\documentclass[journal]{IEEEtran}
\makeatletter
\newcommand{\Spvek}[2][r]{%
  \gdef\@VORNE{1}
  \left(\hskip-\arraycolsep%
    \begin{array}{#1}\vekSp@lten{#2}\end{array}%
  \hskip-\arraycolsep\right)}

\def\vekSp@lten#1{\xvekSp@lten#1;vekL@stLine;}
\def\vekL@stLine{vekL@stLine}
\def\xvekSp@lten#1;{\def\temp{#1}%
  \ifx\temp\vekL@stLine
  \else
    \ifnum\@VORNE=1\gdef\@VORNE{0}
    \else\@arraycr\fi%
    #1%
    \expandafter\xvekSp@lten
  \fi}
\makeatother
\ifCLASSINFOpdf
\usepackage{graphicx}
\usepackage{graphicx,charter, latexsym}

\usepackage{amssymb, amsmath,graphicx,charter, latexsym}
\usepackage{gensymb}
\usepackage{tabularx}

\newtheorem{theorem}{Theorem}

\newtheorem{example}{Example}
\usepackage{multirow,url,float,epstopdf,cite}
\usepackage{verbatim}
\newtheorem{assumption}{Assumption}
\usepackage{graphicx}
\usepackage{wrapfig} 
\usepackage{caption}
\usepackage{subcaption}
\usepackage{algorithmic}
\usepackage{algorithm}
\usepackage{tikz}
\usepackage{pifont}
\usepackage{xcolor}
\usepackage{caption}
\usetikzlibrary{patterns,calc}
\usetikzlibrary{shapes.geometric, arrows}
\usepackage[overlay,absolute]{textpos}
\newcommand{\euler}{e}

\tikzstyle{block}=[draw opacity=0.7,line width=1.4cm]
\tikzstyle{ag} = [circle, radius =3cm, text centered, draw=black]
\tikzstyle{startstop} = [rectangle, rounded corners, minimum width=3cm, minimum height=1cm,text centered, draw=black, fill=red!30]
\tikzstyle{io} = [trapezium, trapezium left angle=70, trapezium right angle=110, minimum width=1cm, minimum height=1cm, text centered, draw=black, fill=blue!30]
\tikzstyle{process} = [rectangle, draw,fill=orange!30, text width = 20em, text centered, rounded corners, minimum height=4em, minimum width=1cm]
\tikzstyle{upd} = [rectangle, draw,fill=orange!30, text width = 5em, text centered, rounded corners, minimum height=4em, minimum width=1cm]
\tikzstyle{decision} = [diamond, minimum width=3cm, minimum height=1cm, text centered, draw=black, fill=green!30]
\tikzstyle{arrow} = [thick,->,>=stealth]
\tikzset{
  solid node/.style={circle,draw,inner sep=1.2,fill=black},
  hollow node/.style={circle,draw,inner sep=1.2},
}

\DeclareTextFontCommand{\emph}{\em}

\IEEEoverridecommandlockouts

\else
\fi
\begin{document}
%
\title{Decentralized Control via Dynamic Stochastic Prices: The Independent System Operator Problem}
%
%
%

\author{Rahul Singh,~\IEEEmembership{Student Member,~IEEE,}
        P.~R.~Kumar,~\IEEEmembership{Fellow,~IEEE,}
        and~Le Xie,~\IEEEmembership{Senior Member,~IEEE}
\thanks{R. Singh is at 32-D716, LIDS, MIT, Cambridge, MA 02139; P. R. Kumar and
Le Xie are at 
Dept. of ECE, Texas A\&M Univ., 3259 TAMU, College Station, TX 77843-3259).
{\tt\small rsingh12@mit.edu, {prk,le.xie}@tamu.edu.}
}
\thanks{Preferred address for correspondence: P. R. Kumar, Dept. of ECE, Texas A\&M Univ., 3259 TAMU, College Station, TX 77843-3259.}
\thanks{This material is based upon work partially supported by NSF under
Contract Nos. ECCS-1546682, NSF Science \& Technology Center Grant CCF-
0939370, NSF ECCS-1150944 and DGE-1303378.}}

%
%

\markboth{}%
{Shell \MakeLowercase{\textit{et al.}}: Decentralized Control via Dynamic Stochastic Prices: The Independent System Operator Problem}
%



\maketitle

\begin{abstract}
A smart grid connects several agents who may be electricity consumers/producers, 
such as wind/solar/storage farms, fossil-fuel plants, industrial/commercial loads, or 
load-serving aggregators, all modeled as stochastic dynamical systems. In each time period, each consumes/supplies some electrical energy. Each agent's utility is either the benefit accrued from its consumption, or the negative of its generation cost. 
There may also be externalities modeled as negative utilities. The sum of all these utilities, called the social welfare, is the total benefit accrued from all consumption minus the total cost of generation and externalities. The Independent System Operator is charged with maximizing the social welfare subject to total generation equalling consumption in each time period, but without the agents revealing their system states, dynamic models or utility functions. It has to announce prices after interacting with agents via bid-price interactions where agents respond with their optimal generation/consumption. 

If agents observe and know the laws of uncertainties affecting other agents, then there is an iterative price-bid interaction that leads to the global maximum value of social welfare attainable if agents had pooled their information. 

In the important case where agents are LQG systems, the bid-price iteration is dramatically simple and tractable, exchanging
only time-vectors of future prices and consumptions/generations at each time step. Agents need not know of the existence of other agents. State-dependent bidding/pricing is not needed. If the DC Power Flow Equations are incorporated
it yields the optimal stochastic dynamic locational marginal prices.

Thereby a solution is proposed for a potentially economically important decentralized stochastic control problem. 
The results may be of broader interest in general equilibrium theory of economics for stochastic dynamic agents.
\end{abstract}

\begin{IEEEkeywords}
Decentralized Stochastic Control, Social Welfare, General Equilibrium Theory,
Demand Response, 
Renewable Energy, Power Systems, Independent System Operator, Energy Market.
\end{IEEEkeywords}

%
\IEEEpeerreviewmaketitle

%
%
%
%

\section{Introduction} \label{sec-intro}
\IEEEPARstart{I}n the electricity grid, the power generated should be equal to the power consumed at all times, neglecting line losses. Unlike other commodities, electricity cannot be stored in the grid. The task of ensuring that generation is balanced with consumption, and in the most economical way, is entrusted to the Independent System Operator (ISO) in deregulated electricity markets
\cite{wu2005power}. 

In the era with fossil fuel as the dominant source of electricity it was possible to adjust generation to meet demand.   In the future, 
as more energy from uncertain and dynamically varying renewables such as wind or solar is used, it is demand that may need to be adjusted continually to balance generation. This strategy is called ``demand response." An example of an adjustable load is an inertial thermal load such as a home with an air conditioner that can be turned off for a while, while still maintaining comfort within the stipulated band of temperatures. New business models are emerging for intermediaries such as retail power service providers, also called ``aggregators" or ``load-serving entities," to sign up a large collection of such customers and undertake their demand response opportunistically, in response to shortages or excess of renewable power that reflect themselves in higher or lower prices, respectively. Large commercial enterprises and industrial loads also will  similarly adjust and optimize their energy usage and cost in-house. Therefore both demand and supply will generally be dynamic and uncertain due to external factors such as uncertain supply and ambient temperature interacting with load requirements.  

The problem we address is how the ISO can perform its task in the new scenario where loads and generators are stochastic dynamic systems. The primary mechanism for coordinating all entities is by time-varying stochastic prices. However, being stochastic dynamic systems, entities will need to know the probability distribution of future prices to plan their optimal consumption/generation over time. But future prices depend on all future uncertainties affecting any of the entity. Uncertainty in wind in a certain locale may affect a wind farm, cloud cover may affect a solar farm, a broken turbine blade may affect a gas turbine, low customer traffic may affect a commercial entity, or high ambient temperature may affect a group of homes, and each of them may globally impact prices everywhere at all future times. However, an entity is generally unaware of uncertainties affecting other entities or how they will respond. So how can they optimally plan their generation/consumption in the face of dynamic uncertainty and lack of knowledge of each other?

The role played by price in coordinating agents has been explored in general equilibrium theory, initiated by Walras~\cite{walras2}. In their breakthrough work, Arrow and Debreu~\cite{arrow,debreu,Debreu1954} showed that a correct choice of prices for commodities ensures, under a quasi-concavity assumption on utility functions, that a system of individual entities, where each optimizes its own response given prices, results in a systemwide Pareto optimal solution where no entity can benefit without another losing. Subsequently, Arrow, Block and Hurwicz  \cite{arrowstable} showed that the prices can be discovered by Walrasian tatonnement~\cite{walras2} under appropriate conditions such as gross substitutes. Their theory extends to allow for uncertainty by simply considering each good under a different random state of nature as a different good, as shown by Arrow~\cite{arrowsecure}. Subsequently, Radner~\cite{radneruncertain} has shown the existence of prices corresponding to an equilibrium even if different agents have different random observations.

The idea of employing prices to perform this task in the electrical power domain was introduced in 
seminal papers by Caramanis, Bohn and Schweppe \cite{caramanis1982optimal}
and Bohn, Caramanis and Schweppe \cite{bohn1984optimal}.
Hogan \cite{Hogan} further elaborated the detailed implementation of a locational marginal price-based electricity market operation. 

Fundamentally based on a static dispatch with no uncertainty, today's electricity market design and corresponding price signal are simply not designed for achieving social welfare optimality for
dynamic generators and loads.
The current market mechanism requires participants to make decoupled bids for separate time intervals. In the day-ahead market, a generator has to bid a price-generation curve for the 8am-9am slot,
another separate curve for the 9am-10am slot, and so on, for each hour of the next day. However, generators have ramping constraints, such as 50 MW/hour, which give rise to inter-temporal constraints between
different time slots. These are typically handled by ad hoc out-of-market (OOM) merit order measures
\cite{ERCOT-Anciliary}.The bidding procedure fundamentally does not allow a generator to bid a
time function even though that is critical to its operation.
Similarly, in the real-time market, the bidding process does not allow a participant to optimize with respect to stochastic process models of uncertain resources such as wind.
There have been many studies on the potential problems associated with this market design, such as unnecessarily price volatility \cite{Roozbehani}, network externalities \cite{Chao}, and  lack
of investment signals \cite{huneault1999review}. While in conventional systems the deterministic and static approach to approximating the underlying dynamic and stochastic power system may be practically appealing without much loss of
optimality, emerging resources such as demand response and intermittent renewables render such approximation invalid \cite{Xie2011windintegration}.
No previous work achieves social optimality of the entire collection
of all stochastic dynamic systems. Providing a theoretical foundation for achieving this fundamental
goal is the target of this paper.

We address the ISO problem where each agent is an individual stochastic dynamic agent whose very nature -- its dynamic model, uncertainties affecting it, and its utility function -- are not necessarily disclosed to others. Our goal is to attain a global maximum of the social welfare. 

The key issue here is not existence of a solution, since here that is simply the maximizer of social welfare, but how to arrive at it,
and realize it, in a distributed way.

There are several interesting aspects to the problem faced by the ISO.
We seek a global optimum of the total social welfare, not just an equilibrium. To see the difference, one can consider the work of Radner \cite{radneruncertain} that is closest to ours in its allowance of different observations for different agents. In that theory, the actions of agents are constant over their information -- if an agent does not know the states of other agents, then its action does not change unless its own observed state changes. However, a globally optimal solution will require coordination of the actions of all entities so as to be responsive to each others' states. Thus, the price stochastic process will need to provide this additional coordinating information. The issue examined in this paper is how the ISO is to determine prices and ensure such coordination. 

Importantly, we will fundamentally exploit the very fact that uncertain events unfold over time in a dynamic system, to design dynamic interactive strategies for coordination. This is in contrast to Arrow's approach~\cite{arrowsecure} where the problem with uncertainty is reduced to a problem without uncertainty when the descriptions of the uncertainty states of all the agents and their utility functions are known a priori. The equilibrium prices corresponding to each uncertainty state can then be computed at the very outset itself. Such an approach therefore considers the problem in ``normal" form, where the entire dynamic system is simply formulated as a ``static" system where each agent chooses its strategy as a function of prices at  states. This observation is also made by Smale~\cite{smale}. 

The problem is also interesting from the viewpoint of decentralized stochastic control. Since we seek to maximize social welfare, it is a problem in team theory. However, since Witsenhausen~\cite{wsh}  it is known that if agents are unaware of each other's actions but influence each other's observations, then the problem is generally intractable, even in linear quadratic Gaussian (LQG) systems. Unawareness of other's actions is the norm in any distributed stochastic system such as the ISO problem.
Nevertheless, in the ISO problem, we show that a system consisting of a collection
of LQG systems has an elegant and tractable solution.
The tatonnement process for obtaining the global optimum is remarkably simple, even if all entities have private uncertainties. There is no need for agents to even share models of their systems, their uncertainties, the probability distributions of their uncertainties, or their utility functions. Thus, the ISO can optimally coordinate a set of distributed LQG systems very simply without knowing any of their details. This is potentially important, since LQG models are widely used in power systems, where systems are often approximable as linear systems, noises as Gaussian, and costs as quadratic in states and actions.

Importantly, the above approach extends to any number of linear constraints, besides balancing generation and consumption. An important task of the ISO problem is to ensure delivery of required power flows over a congested transmission network. A commonly used of model the transmission network, is through (somewhat misleadingly labeled) ``DC power flow" equations, where differences in bus phase angles determine the power flows along the lines~\cite{bergen}. Their popularity derives from the fact that the resulting equations are linear. Hence the LQG model extends to include the transmission network and provides a very simple solution for the ISO to obtain dynamic \emph{stochastic dynamic} locational marginal prices that attain the global maximum of the social welfare.

The paper is organized as follows.
Section \ref{sec-relatedworks} surveys related work.
Section \ref{sec-dynamic} describes the broad context,
Section \ref{sec-model} the system of agents,
formulates the ISO problem, and
describes the fundamental challenges.
We then progressively build up to more complex systems.
Beginning with static deterministic systems in Section \ref{sec-static}, we
show how the ISO can determine both the optimal price and optimal allocations of consumptions/generations
to agents through a bid-price process that corresponds to a subgradient iteration with subsequent averaging.
Then we consider
deterministic dynamic systems in Section \ref{sec-deterministic} 
and show how prices and allocations as a function of time can be determined through
the same bid-price iteration. 
Section \ref{sec-ibs} describes iterative bid-price schemes used subsequently in the stochastic dynamic context.
Then we turn to the stochastic problem
in Section \ref{sec-common} where all agents are subject to a common uncertainty
and show that by viewing the system as a ``tree" the results can be extended from the deterministic case.
We next show in Section
\ref{sec-private} that this approach can be extended to systems where entities have private uncertainties,
by viewing the system in extended form, and iterating at each time between the ISO and the agents for price 
and allocation discovery.  Knowledge by the agents of the probability laws of each other's uncertainties
is required, though not of their dynamic models, utilities, or the semantics of the uncertainties
since labels of uncertainties can be non-informatively chosen. The difficult issue is complexity, caused by the 
exponentially exploding joint state space of the uncertainties. However, in Section \ref{sec-lqg}, we show that
the complexity disappears for distributed LQG systems, leading to a simple and implementable solution.
The ISO simply discovers and announces time-varying but \emph{not} state-dependent prices for future epochs,
and revises them at each time step, reminiscent of model predictive control. 
We show in Section \ref{sec-powerflow} that this result can be extended to include any linear constraints, e.g., the widely used DC-power flow equations.
Section \ref{sec-simu} presents 
the results of illustrative simulations, concluding in Section \ref{sec-concluding}.

\section{Related Works} \label{sec-relatedworks}
No similar results appear to be known for general decentralized stochastic control. Team problems have been extensively studied,
e.g.,~\cite{radnor,sardar,van}, but those formulations do not apply here since agents need to know the system dynamics of other agents. Even when the models are known, there are still considerable difficulties in decentralized stochastic control. When agents do not share observations, severe complexity arises, even in LQG systems, as shown by Witsenhausen's counterexample of a two stage problem~\cite{wsh}. The roles of observation, signaling~\cite{sardar}, and the trade-off between communication and control are evident from Witsenhausen's counterexample~\cite{wsh}. Teneketzis~\cite{demos} considers decentralized stochastic control under the restrictive assumption that the interaction between agents is ``weak". There are some recent structural results~\cite{nayar}, and results regarding sufficient statistics~\cite{jeff} under these assumptions.

From the economics side, this work is an extension of general equilibrium theory~\cite{arrow1}. To the authors' knowledge there does not appear to be any similar result for coordinating multiple LQG systems or the efficiency of the simplified signaling. While the name may appear to be related to the issues studies here, Dynamic General Stochastic Equilibrium theory pioneered in \cite{kydland1982time} addresses issues in macroeconomics, and is not relevant for the problems of interest here.

Viewed from the power system end, there have been many efforts since the deregulation of the electricity sector on a market-based framework to clear the system. Today's locational marginal price-based nodal market design is based on seminal work in \cite{bohn1984optimal,Hogan}. This has been followed up
by a large body of literature focusing on designing an efficient transmission pricing mechanism in support of an efficient market \cite{wu1996folk} \cite{kirschen1997contributions}. From the system operators' perspective, the naive belief that deregulation of electricity industry would simply work was critically re-assessed following 
the Enron crisis and lack of long-term investment \cite{lave2004rethinking} \cite{hogan2003transmission}. From a market participant's perspective, there has been pioneering work on game theoretic approaches to modeling the market power issues in the electricity market \cite{chuang2001game} \cite{baldick2004theory}. With increasing penetration of stochastic resources, there have been efforts at designing a market bidding mechanism that achieves the social welfare optimum. Ilic et al. \cite{Ilic2011framework} have proposed a two-layered approach that internalizes individual constraints of market participants while allowing the ISO to manage the spatial complexity. References~\cite{carpenter,carpenter1} contain some heuristic approaches. Reference~\cite{wets} applies progressive hedging to deal with uncertainties on the production side; however the solution is centralized, and does not provide any theoretical guarantees.
Reference \cite{rajagopal2013risk} studies how the ISO should dispatch, i.e., purchase energy and call options in different markets,
under forecast errors about future loads and renewable generation, when future  decisions can mitigate current errors.

However, there has been no analytical framework that precisely leads to the social welfare optimality with dynamic, stochastic inputs from market participants. The major challenge addressed is how to elicit optimal demand response in such cases without
generators/loads revealing the details of their dynamic models to the ISO. 


\section{The ISO Problem of Coordinating Dynamical and Uncertain Generators, Loads and Prosumers} \label{sec-dynamic}
Generators such as wind farms, photo-voltaic farms, hydro, coal or gas turbines, need to be modeled as dynamic stochastic systems. Likewise, loads such as aggregators, commercial or industrial establishments, also need to be so modeled since their demand may be governed by a dynamic system, and random. Since environmental variables such as temperature are involved, and since human beings in the loop respond to economic incentives/prices, their response may also be uncertain. Hence loads generally will also need to modeled as stochastic dynamic systems. Storage services such as battery farms or pumped hydro can also be modeled as dynamic systems where the state is the amount of energy stored. Dynamic models can also be used to model ``prosumers" such as homes with solar panels, which can switch between consumption/generation.

The utility of a generator is the negative of its cost of generation. The utility of a load is the ``benefit" that the load accrues from the consumed power. There may also be externalities, such as pollution, that can optionally be modeled as a cost, i.e., negative utility. The total of all the utilities, called social welfare, is therefore the benefit of the power consumed minus the cost of generating it and the cost of externalities, and the goal is to operate the overall system so as to maximize it. 

An agent's utility is measured not statically, but over a time interval of interest, since agents are dynamic systems. A large commercial load may accrue utility over a period of time, by maintaining the temperature within a band by switching air conditioners off and on taking time lags into account, if demand response strategies are in place. Similarly, a storage service may buy and store energy at off-peak times and sell it at peak times, again accruing value only over a time interval. Generators too may accrue utility over time by ramping up generation.

An agent's utility is affected by stochastic uncertainties of other agents, since coal shortage at a generator may affect a distant load.
Each agent seeks to
maximize the expected value of its own utility function, with expectation taken over all uncertainties affecting all agents.

There are important and severe constraints on the information disclosed to the ISO. It does not know the states, dynamic models, or utility functions of individual agents. Whether loads or generators, 
individual agents may be averse to disclosing information to others for competitive reasons
or to ensure privacy. A load-serving entity may not inform the ISO of the states of its loads, e.g., the temperatures of every one of its large collection of customer's homes. A solar farm may not inform the ISO of the extent of its cloud cover. More fundamentally, the agents may not even be willing to share their individual dynamic system models with each other. The ISO may not be informed of the stochastic model of wind at a wind farm, or the detailed model of the dynamics of a coal plant. Similarly, agents may not share their individual utility functions. A load-serving entity may not be willing to disclose its contracts with its customers or its cost of operations. Many of these entities compete with
each other and are sensitive to sharing their information with competitors. 
Recently, privacy of loads has also emerged as a major concern.\footnote{Interestingly, privacy is nowhere mentioned in the seminal paper~\cite{caramanis1982optimal} that introduced spot prices, indicative of how new issues arise.}  Demand response can entail violation of privacy since it has been shown that having access to a home's real power consumption allows one to deduce the number and behavior of its occupants~\cite{molina2010private}.

Even if all agents were willing to share all their information with the ISO, it would be such an intractably large amount of information, amounting to a complete state of the world, that the ISO would not be able to handle it with acceptable complexity and delay anyway.

The ISO is nevertheless charged with maximizing the sum of the utilities all the agents, i.e., the social welfare. It plays the role of a mediator. It needs to determine how much power each agent should generate/consume in
each time period over the time horizon of interest. It needs to
allocate the required power generation to the generators over time in the most economical fashion.
It also needs to provide the optimal amount of power over time to each load to optimize its utility.
Demand and supply are intertwined, since demand is uncertain and determined by the cost of generation, 
while generation is also uncertain and incentivized by the price that consumers are willing to pay.
All the agents are stochastic dynamic systems with their own utility functions over time. 
The ISO needs to do all this in the face of all the stochastic uncertainties affecting the agents over the time horizon. 

The ISO would like to achieve the above through economic mechanisms, by determining prices, and by agents
responding with their own selfish utility maximizations as in general equilibrium theory~\cite{arrow}. The complication is that the agents are evolving stochastically in time. Therefore the prices cannot simply be announced once and for all at the beginning of the time interval of interest. For example, a future high ambient temperature can lead to very high demand that will need high prices to incentivize extra generation or reduce deferrable demand. The prices will need to vary stochastically in response to private uncertainties affecting the agents.
The price stochastic process should carry all the information that is necessary for all entities in the overall system to coordinate in a globally optimal way, since each affects others.

The fundamental question examined in this paper is whether and how this optimality can be attained given stochastic dynamical system models for the agents, and what form the mediation process or tatonnement~\cite{walras1} takes. Our contribution is to show that there are iterative interaction processes under which the ISO can indeed perform this task. We address the complexity of this task under several scenarios. The complexity is very high in the general case. However, in the case where the agents can be modeled as linear Gaussian stochastic systems and the cost functions are quadratic, we show that a much simpler scheme yields the systemwide global optimum. This scheme extends to encompass other linear constraints, e.g., those modeling the electrical network, thus providing a more comprehensive solution that takes into consideration the power flows over the network.

Beyond balancing generation and consumption, there are at least two additional problems that the ISO faces.
It needs to ensure no line's capacity is exceeded in the electrical transmission network, so as
to prevent overheating. This requires ensuring that in the solution of the power flow equations at
the obtained generations, the current carried over every line does not exceed its capacity~\cite{bergen}.

ISO's also need to ensure reliability. If contingencies occur, such as generator tripping or line-to-ground fault, then the system's electrical state should converge to an acceptable equilibrium point \cite{gan2000stability}. ISO's verify that the solution has this property for all single event contingencies, reformulating the problem if any violations are observed.  Considering multiple simultaneous contingencies is computationally demanding and not the norm \cite{wang2013risk}.

\section{System Model and ISO Problem}\label{systemmodel} \label{sec-model}
We consider a smart grid consisting of $M$ agents, each of which may act as a producer, consumer or both, i.e., a prosumer, evolving
over a time interval $t=0,1,\ldots,T-1$.  
The time horizon $T$ could be $96$, which would correspond to one day of
$15$ minute slots in the real-time market.
There is however considerable flexibility to model other scenarios. 
One can model the risk-limited dispatch of \cite{rajagopal2013risk} 
where purchases of forward energy are made for
blocks of time, with blocks getting shorter as operations approach real time.
In that case the times $t=0,1, \ldots $ can correspond to the 24 hour ahead, all the one-hour ahead, and
all the 15 minute ahead, times
at which decisions are made by agents.
The states of the agents (below) can keep track of their past purchases, temperature forecasts, etc,
so that noises can be regarded as changes to past forecasts, allowing considerable generality.
\\
%
%
${\mathbf{Randomness}}$ is modeled through a probability space $\left(\Omega,\mathcal{F},\mathbb{P} \right)$. The ``state of the world" $\omega \in \Omega$ captures a multitude of random phenomena spread out temporally and spatially, for example, unpredicted weather (the wind-speed), unexpected events such as coal shortage, or a damaged wind-turbine, etc. \\
{{\bf Common and Private Uncertainties.}} The randomness $\omega$ affects an agent $i$ through the stochastic processes $N_i(\omega,t)$ and $N_c(\omega,t), 0 \leq t \leq T-1$. $N_c(t)$ is a ``common" uncertainty that affects and is known causally by all agents, e.g., the weather of a city. $N_i(t)$ is a ``private" uncertainty that is specific to agent $i$, and known causally only to agent $i$. As we will see, this decomposition of uncertainties clarifies the task of constructing interaction schemes between the agents and ISO. \\
$\mathbf{Agents}$ are modeled as  stochastic dynamical systems.
The state $X_i(t)$ of agent $i$ at $t$ is known to it, and evolves as
$
X_i(t+1) = f_i(X_i(t),U_i(t),N_i(t),N_c(t),t)
$, where $f_i$ describes the dynamics of the agent $i$. 
The initial condition $X_i(0)$ can be random.
The common case is that $U_i(t)$ is a scalar that denotes the amount of electricity consumed (negative if supplied) from the grid by agent $i$ at time $t$, but we will allow $U_i$ to be a vector of several commodities being produced and consumed. \\
${\mathbf{Consumption/Generation \ Constraints}}$.
Let $N_i^t := (N_i(0), N_i(1), \ldots , N_i(t)$ denote the past of $N_i$, and similarly define $N_c^t$.
Agent $i$'s choice has to satisfy the local capacity constraints
$F_i(N_i^t,N_c^t,t)U_i(t) \leq g_i(N_i^t,N_c^t,t)) + \sum_{s=0}^{t-1} C_{i}(N_i^t,N_c^t,s,t)U_i(s)$ and
$h_i(N_i^t,N_c^t,t,U_i(t)) \leq 0$ for each $t$. 
The affineness of the former constraints in past $U_i$'s allows for ramping,
the dependence on $t$ allows for seasonality, and the dependence on $N_i,N_c$ allows
random effects on capacity. 
 \\
${\mathbf{Observations}}$ available to an agent $i$ until time $t$ include the realizations of its system state $X_i(s)$, common noise $N_c(s)$, and its private noise $N_i(s)$ for $0 \leq s\leq t$.\\
The {\bf{One-step Cost Function}} of an agent $i$, $1 \leq i \leq M$, denoted $c_i(x_i,u_i,t)$ (or its negative, a one-step utility function $-c_i(x_i,u_i,t)$), is a function of its state and action, in period $t$. For producers, this could be the cost of labor, coal, etc.. For consumers, this could represent the cost incurred due to the high temperature of house/business facility, or due to a delay in performing a task resulting from inadequate purchase of electricity, or the negative of some benefit
of the electricity usage. \\
{\bf{Externalities}}, e.g., pollution, with one-step cost
$\sum_{i=1}^M e_i(u_i,t) $, say cost of mitigation, can be considered. 
By allowing $e_i(u_i,t)$ to be positive/negative ISO imposed levies, cross-subsidies can be addressed.
For linear levies, $e_i u_i(t)$, ISO budget balance, $\sum_{i=1}^M e_i u_i(t) = 0$, can be enforced as 
shown below for energy balance. \\
{\bf{Energy Balance}} should be maintained in each period, i.e., $\sum_{i=1}^M U_i(t) = 0$ for all $t=0, 1, \ldots , T-1$. 
We allow {{\bf general linear vector constraints}}: $\sum_{i=1}^ME_i(t)U_i(t)=d(t)$.  \\
{\bf{Total System Operating Cost}}, or its negative, the {\bf{Social Welfare}}, is  the sum of
the expected value of the 
finite horizon total of the one-step costs incurred by all the agents plus externalities,
$
\mathbb{E} \sum_{t=0}^{T-1} \sum_{i=1}^{M+1}[c_i \left(X_i(t),U_i(t),t \right)+e_i \left(U_i(t),t \right)]$.
It is the total electricity generation cost plus the cost of externalities minus the utility provided to the consumers. The expectation above is taken with respect to the combined uncertainty or ``noise"€ process $ N(t) := \left(N_c(t), N_1(t), N_2(t), \ldots, N_M(t)\right)$ for $t=0,1, \ldots, T-1$, consisting of all the private uncertainties and the common uncertainties, as well as the random initial conditions of all the $M$ agents. \\
The {\bf{Power Flow Equations}} are algebraic equations based on Kirchoff's laws that have to be satisfied by the electrical variables, voltage and current magnitudes and phase angles. They impose constraints on $\{U_i(t), 1 \leq i \leq M \}$. 
\\
The {\bf{Independent System Operator (ISO)}}'s task is to maximize the social welfare.
It solicits electricity purchase/sale bids from the agents in each time slot $t=0,1,\ldots,T-1$. Our model allows for agents and the ISO
to iterate on the bids. Once the price iterations have converged, the ISO declares the market clearing prices, and the 
electrical energy to be consumed/generated by the agents, at the declared prices.   \\
The {\bf{Bidding Schemes}} allow the ISO and agents to reach a solution for prices, generation and consumption.
Depending on the assumptions made about the system model, there will be different bidding schemes. 
An example is the following.
Consider time $s$.
The ISO announces a price sequence for current and future times, $s \leq t \leq T-1$, to all agents.
Agent $i$ bids, as a function of its past information, the amount of electricity it is willing to purchase/generate
at the current and future times $s \leq t \leq T-1$,
at the prices indicated by the ISO.
After collecting the bids, the ISO updates the price sequence.
An iteration of \textit{price updates} followed by \textit{bid updates}, continues till the prices and the bids converge,
and then the ISO announces the allocations of generations/consumptions to agents for the current period $s$. 
This entire process can be repeated in each discrete time slot $s$ in real-time. 
\\
 {\bf{Goal of Social Welfare Maximization}}: The goal is to maximize the negative of total system cost, i.e., social welfare.  Let $\mathcal{F}_t$ be the $\sigma$-algebra
generated by all the noises upto time $t$, as well as initial conditions. Now we come to the stringent goal of this paper. \emph{We would like to attain the same maximum value of the social welfare as could be attained over the class of all control laws where $U(t) := \left(U_1(t), U_2(t), \ldots , U_M(t)\right)$ is adapted to $\mathcal{F}_t$}. 

This is an 
economically important point in that even though the agents do not all act in a centralized way and even though they do not all have access to all the observations and initial conditions of each other, we would like them to collectively attain the same optimal value of social welfare by acting in a distributed way, with each agent only using its own causal observations together with the price information provided by the ISO. In fact they do not even know each other's dynamic models or cost functions, taking this problem outside of usual stochastic control/game theory. 

The resulting {\bf{ISO Problem}} is:
\begin{align}
& \min \mathbb{E} \sum_{t=0}^{T-1} \sum_{i=1}^{M} [c_i \left(X_i(t),U_i(t),t\right) + e_i \left(U_i(t),t\right)]  \label{p0} \\
&\mbox{such that } X_i(t+1) = f_i(X_i(t),U_i(t),N_i(t),N_c(t), t); \notag \\
& \mbox{with capacity constraints } h_i(N_i^t,N_c^t,t,U_i(t)) \leq 0, \label{constraints-at-time-t-1} \\
& F_i(N_i^t,N_c^t,t)U_i(t) \leq g_i(N_i^t,N_c^t,t) \notag \\
&\quad \quad \quad \quad \quad \quad \quad \quad \quad + \sum_{s=0}^{t-1} C_{i}(N_i^t,N_c^t,s,t)U_i(s); \label{constraints-at-time-t-2}\\
&\sum_{i=1}^ME_i(t)U_i(t)=d(t) \mbox{ for }1 \leq i \leq M, 0 \leq t \leq T-1.  \label{constraint-on-balancing}
\end{align} 
The expectation above is taken with respect to the combined uncertainty or ``noise" process $N(t):=\left(N_c(t),N_1(t),N_2(t),\ldots,N_M(t),N_M(t)\right), t=0,1,\ldots,T-1$, as well as the random initial conditions $\left(X_1(0), X_2(0), \ldots , X_M(0)\right)$. The actions $U_i(t)$ are to be taken on the basis of the past information available to agent $i$ at time $t$, which includes the past history of its own observations of its system's state, common noise and private noise, as well as the common price information provided to all agents by the ISO. The bidding process to be studied below will determine the prices announced by the ISO to all the agents.

The central issue is the following: \emph{How should the ISO determine pricing and allocations to dynamic stochastic agents so that the overall system is as optimal as it could be through centralized control, even though agents do not know each other's dynamic models or cost functions?}

\subsection{Fundamental Challenges} \label{sec-fi}
The ISO Problem
poses several challenges. It is a multi-agent problem where stochastic dynamic agents with differing objectives,
ignorance of each other's systems or objectives, and separate observations, are constrained in their joint actions; yet the goal is to ensure that they function as a team and jointly maximize social welfare. 
\subsubsection{Constraint on joint actions}\label{fi:id}
The problem cannot be solved by considering each agent separately because energy balance (\ref{constraint-on-balancing}) constrains their joint actions.
\subsubsection{Privacy constraints}\label{privacy}
The agents do not disclose their system dynamics
functions $f_i$ to other agents or the ISO.
In fact, the agents do not even know how many agents are present.
\subsubsection{Non-classical information structure}\label{fi:dc}
Even if all dynamics and probability distributions of uncertainties were known to all, the ISO Problem would still lie at the core of decentralized stochastic control with a non-classical information structure
\cite{wsh,wsh2,sardar,van,demos,carpenter,nayar} since each agent has separate observations from others.
Even if privacy were not an issue, sharing all observations amongst all agents requires huge communication and processing overhead, etc., and may be impossible in practice. 
\subsubsection{Conflicting objectives} The objectives of the agents are not all aligned and may have conflicts.
\subsubsection{Signaling}
In decentralized stochastic control~\cite{sandelv,nayar,van},
controllers may be able to signal some private information to other agents
over a ``channel" which may even be the physical
plant itself. 
``Prices" can play the role of a channel, with the bidding scheme functioning as encoder-decoder. Essentially the ISO needs to construct a ``price" sufficient statistic for the problem~\cite{striebel,blackwell}. 

The question we address is: \emph{Can $M$ independent systems be driven to an overall optimal operation through the ISO}? We will show that there exist ``iterative bidding schemes" (IBS) which yield the same performance as that of the optimal centralized controller.

\section{Static Deterministic Systems} \label{sec-static}
We begin with the problem where all generators and consumers are static and deterministic. 
The ISO 
has to allocate generations and consumptions so that the social welfare, the total benefit accrued from consumption
minus the total cost of generation is maximized.

This can be formulated as a problem of minimizing the total cost $J(u) = \sum_{i=1}^M [c_i(u_i) + e_i(u_i)] $ of $M$ agents
and the externalities, 
where, if agent $i$ is a generator producing $-u_i$ (negative by convention for generation), the cost of generation is $c_i(u_i)$,
while if it is a consumer consuming $u_i$ the utility it obtains from consumption is $-c_i(u_i)$, and $e_i(u_i)$ is the externality
associated with generation/consumption. 
Each generator/consumer $i$ may also be subject to linear/nonlinear vector inequality constraints $F_i u_i \leq g_i$ and $h_i(u_i) \leq 0$. (When $u_i$ is a scalar, this will reduce to either a semi-infinite interval or interval constraint on $u_i$ under convexity of $h_i$ below).

This entails solving the following optimization problem:
\begin{align} 
&\min_{u_1, \ldots , u_M} \sum_{i=1}^M [c_i(u_i) + e_i(u_i)], \label{simple-primal} \\
&\mbox{subject to: } F_i u_i \leq g_i, h_i( u_i) \leq 0  \quad \mbox{ for } 1 \leq i \leq M, \label{individual-constraints} \\
& \mbox{and } \sum_{i=1}^M E_i u_i = d. \label{balancing_constraint}
\end{align}

Dualizing only the constraint (\ref{balancing_constraint}), and denoting $u := (u_1, u_2, \ldots, u_M)$, yields, respectively, the Lagrangian,
Dual Function, and optimal reward of the Dual Problem:
\begin{align}  
&\mathcal{L}\left(u, \lambda \right) := \sum_{i=1}^{M} [ c_i(u_i) +  + e_i(u_i) + \lambda^T E_i u_i  ] -\lambda^Td, \notag \\
&D(\lambda) := \min_{ \{u:  F_i u_i = g_i, h_i(u_i) \leq 0 \  \forall i \} }  \mathcal{L}\left(u, \lambda \right), \notag \\
&J^\star := \max_{ \lambda } D(\lambda) = D(\lambda^\star) \label{Simple_Dual} .
\end{align}

\begin{assumption}[Assumption for deterministic case]\label{finite-compact-slater}
(i) $c_i( \cdot)$, $e_i( \cdot )$, $h_i( \cdot )$ are convex, $\{u_i:  F_i u_i \leq g_i, h_i(u_i) \leq 0 \}$ is compact,
and (\ref{simple-primal},\ref{individual-constraints},\ref{balancing_constraint}) has an optimal solution. \\
(ii) Slater's Condition: There exists a feasible $\bar{u}_i$ satisfying $h_i( \bar{u}_i) < 0$
in $\mbox{RelInt(Dom}(c_i)) \cap \mbox{RelInt(Dom}(e_i))$. 
\end{assumption}
From (ii), $J^\star$ is
also the optimal cost of the Primal (\ref{simple-primal}).

Since $D(\lambda)$ can be decomposed agent-by-agent as
\begin{align*}  
D(\lambda) = \sum_{i=1}^M \min_{ \{u_i:  \mbox{ \scriptsize{s.t. (\ref{individual-constraints})}} \} }  [ c_i(u_i) + e_i(u_i) + \lambda^T E_i u_i ]-\lambda^Td,
\end{align*}
the ISO can conceivably simply announce the ``optimal price" $\lambda^{\star}$ 
per unit of power as that which attains the max in (\ref{Simple_Dual}), and assess an additional levy $e_i(u_i)$ on agent $i$.
(This levy could be a ``carbon tax" used to mitigate the pollution). 
Each agent $i$ can then respond with either its generation $-u^*_i$ or consumption $u^*_i$
that minimizes its net ``loss" $c_i(u_i) + \lambda^{\star} u_i$ over (\ref{individual-constraints}).
The ISO can finally announce the generation/consumption allocations to the agents.

There are two issues that arise: \\
(i) Since agents do not disclose their cost functions, there needs to be a price discovery process, as in a Walrasian auction~\cite{walras1}. The ISO's price needs to be reduced/increased according to whether the agents' response
results in
excess total generation/consumption). We consider the following
\emph{iterative price-bid process:}
\begin{align*}  
 \lambda^{k+1} & = \lambda^k + \frac{1}{k}[\sum_{i=1}^M E_i u^k_i - d],  \\
u^{k+1}_i & =  
\scriptsize{
\begin{matrix} 
\mbox{argmin} \\ 
{ \{u_i:  \mbox{ \scriptsize{s.t. (\ref{individual-constraints})}} \} }
\end{matrix} 
}
[c_i(u_i) +e_i(u_i) + (\lambda^{k+1})^T E_i u_i )]. \notag
\end{align*}
This iteration of prices\footnote{The gain $\frac{1}{k}$ can be replaced by $\frac{\alpha}{k^\delta}$ for $\frac{1}{2} < \delta \leq 1$
with $\alpha > 0$ in Sections \ref{sec-static}--\ref{sec-private}.} and bids is a subgradient algorithm that converges to an optimal price
for the Dual under Assumption \ref{finite-compact-slater}
\cite{anstreicher2009two}. \\
(ii) The recovery of optimal generations/consumptions from optimal price is more problematic:
\begin{example}[Counterexample to generation/consumption recovery from optimal price] \label{counterexample}
Consider one generator and one load.
The generator's cost of producing $-u_1$ units of energy is $-\frac{2}{5}u_1$, with $u_1$ restricted to
$[-1, 0]$, and the cost of the externality is $-\frac{1}{10}u_1$. The load's utility from consuming $u_2$ units of energy is $\log(1+u_2)$ with $u_2$
restricted to $[0,2]$, and it has no externality. Energy should be balanced. The social welfare problem is:
\begin{align*}
&\mbox{Min} \quad -\frac{2}{5}u_1 -\frac{1}{10}u_1 - \log (1+u_2) \\
&\mbox{Subject to: }
-1 \leq u_1 \leq 0, 0 \leq u_2 \leq 2, u_1 + u_2 =0.
\end{align*}
The optimal solution is $(u_1^\star,u_2^\star) = (-1, 1)$.

The Dual function of price $\lambda$ is
\begin{align*}
&D(\lambda) = \min_{1 \leq u_1 \leq 0}[ -\frac{1}{2}u_1 +\lambda u_1 ] + [ - \log(1+u_2) ] + \lambda u_2 ].
\end{align*}
The minimizers and minimum, $(-u_1(\lambda),u_2(\lambda), D(\lambda))$, are
$$
=
\begin{cases}
(0,\mbox{Min} \{ \frac{1}{\lambda} - 1, 2 \},1 - \lambda + \log \lambda) \quad \mbox{ if } \lambda < \frac{1}{2},\\
(\mbox{\emph{any} point in } [0,1],\frac{1}{\lambda} - 1, 1 - \lambda + \log \lambda) \quad \mbox{ if } \lambda = \frac{1}{2},\\
(1, \frac{1}{\lambda} - 1, \frac{1}{2} +  \log \lambda) \quad \mbox{ if } \frac{1}{2} < \lambda \leq 1,\\
(1,0, - \lambda +  \frac{1}{2}) \quad \mbox{ if } 1 < \lambda,
\end{cases}
$$
The optimal solution of the Dual is $\lambda^\star = \frac{1}{2}$.

However, when the price $\lambda^\star = \frac{1}{2}$ is announced by the ISO, the generator can bid $-u_1 =0$
since any point in [0,1] is optimal.
The load's bid is $u_2 = 1$, and there will not be balance between generation and consumption.
\hfill $\square$
\end{example}

Therefore one cannot recover the optimal bids from the optimal prices.
However, they can be recovered from the \emph{iterations of the bidding process} under
Assumption \ref{finite-compact-slater}
by taking weighted averages of previous  bids
\cite{gustavsson2015primal}.
Thus the very \emph{process} of iterative bidding is itself important.
\begin{theorem}[Determining optimal bids by generators and loads \cite{gustavsson2015primal}]\label{ergodic-theorem}
Let $\theta \geq 0$.
Consider the averaged bids obtained recursively as follows:
\begin{align}
\bar{u}^k_i = \frac{\sum_{s=1}^{k-1}s^\theta}{\sum_{s=1}^{k}s^\theta} \bar{u}^{k-1}_i + 
\frac{k^\theta}{\sum_{s=1}^{k}s^\theta} u^{k}_i; \quad \bar{u}^0_i = u^0_i  \label{ergodic_method}
\end{align}
Then $\bar{u}^k_i \to u_i^\star$ which is optimal for  (\ref{simple-primal}). \hfill $\square$
\end{theorem}
A larger $\theta$
weights more recent values of the iterates for $u_i$ more heavily, while $\theta=0$ takes a plain average.

\begin{example}[Continued] \label{counterexample-continued}
Choosing $\theta = 2$, one obtains:
\begin{align*}
&\{\lambda_k \}: 0, 0, 1,0.6667, 0.5416,
\ldots \to \frac{1}{2}, \\
&\{u^k \}: \Spvek{-1;0}, \Spvek[l]{0;2}, \Spvek[c]{-0.9412;0.1176}, \Spvek{-0.9898;0.4133} , \\
&\quad \quad \Spvek{0-0.9972;0.7263} , 
\Spvek[l]{-0.9990;0.9526}, 
\ldots \to \Spvek{-1;1}. \hfill \quad \square
\end{align*}
\end{example}

\section{Dynamic deterministic systems} \label{sec-deterministic}
We consider the ISO Problem for deterministic systems:
\begin{align}
& \min \sum_{t=0}^{T-1} \sum_{i=1}^{M}[c_i \left(x_i(t),u_i(t),t\right) + e_i \left(u_i(t),t\right)] \label{deterministic-optimal-control} \\
&\mbox{s.t. } x_i(t+1) = f_i(x_i(t),u_i(t),t), 
\mbox{and (\ref{constraints-at-time-t-1},\ref{constraints-at-time-t-2},\ref{constraint-on-balancing})}. \label{systemconstraint} 
\end{align} 
Since the state variables $x_i(t)$ can be expressed in terms of the inputs $u_i:=\left(u_i(0),u_i(1),\ldots,u_i(T-1)\right)$,
the ISO problem 
can be written as (\ref{simple-primal}).
We assume Assumption \ref{finite-compact-slater}.


The associated Lagrangian and dual function are
\begin{align}  
&\mathcal{L}\left(u, \lambda \right): = \sum_{i=1}^{M} \sum_{t=0}^{T-1} [c_i(x_i(t),u_i(t),t) + e_i(u_i(t),t) \notag \\
& \quad \quad \quad \quad \quad \quad \quad \quad + \lambda (t)^T E_i(t) u_i(t) ] - \sum_{t=0}^{T-1}  \lambda (t)^T d(t), \notag\\
&D(\lambda) := \min_{\scriptsize{\{u: u \mbox{ satisfies } (\ref{constraints-at-time-t-1},\ref{constraints-at-time-t-2})} 
(N_i,N_c \mbox{\scriptsize{ absent) }} \mbox{\scriptsize{for} } \scriptsize{0 \leq t \leq T-1, \forall i\}}}
\mathcal{L}\left(u, \lambda \right), \label{dual-for-deterministic-dynamic}
\end{align}
where $u :=\left(u_1,\ldots,u_M \right)$, $\lambda :=\left(\lambda(0),\ldots,\lambda(T-1)\right)$, and each $x_i(t)$ is regarded as a function of $u_i$ in (\ref{dual-for-deterministic-dynamic}). The Lagrangian decomposes by agents.
Therefore, given $\lambda$, 
each agent $i$ solves its own problem:
\begin{align}
& \mbox{Min} \sum_{t=0}^{T-1} [ c_i(x_i(t),u_i(t),t) + e_i(u_i(t),t) + \lambda(t)^T E_i(t) u_i(t) ] \label{adjoinedcost} \\
&\mbox{subject to } (\ref{systemconstraint}). \notag
\end{align}

The Bid-Price Iteration proceeds as follows.
The ISO announces prices $\lambda^k = \{ \lambda^k(t): 0 \leq t \leq T-1 \}$.
Each agent $i$ responds with an optimal solution $u_i^k:=\left(u_i^k(0),u_i^k(1),\ldots,u_i^k(T-1)\right)$ to (\ref{adjoinedcost}).
Since the subgradient with respect to $\lambda$ of the Dual function $D(\lambda)$ is $\left(\sum_{i=1}^M u_i(0),\sum_{i=1}^M u_i(1),\ldots,\sum_{i=1}^M u_i(T-1)\right)$, the ISO employs the price iteration over $k$, for every $t \in [0, T-1]$:
\begin{align}
&\lambda^{k+1}(t) = \lambda^k(t) + \frac{1}{k} [\sum_{i=1}^M E_i(t) u^k_i(t) - d(t)]. \label{PriceIterate}
\end{align}
The agent bids are averaged by (\ref{ergodic_method}) to give $\bar{u}^k_i(t)$. The ISO announces the allocations $u_i^\star(t) := \lim_{k \to \infty} \bar{u}^k_i(t)$.

\begin{theorem} \label{theorem-for-determinstic-dynamic-case}
Consider the ISO problem (\ref{deterministic-optimal-control},\ref{systemconstraint}) under Assumption \ref{finite-compact-slater}.
Suppose the ISO employs the price iteration (\ref{PriceIterate}), with each agent
$i$ responding with an optimal solution $u_i^k:=\left(u_i^k(0),u_i^k(1),\ldots,u_i^k(T-1)\right)$ to (\ref{adjoinedcost},\ref{systemconstraint}).
Then the prices $\lambda^k$ converge to the optimal prices $\lambda^\star$
for (\ref{dual-for-deterministic-dynamic}).
The ISO's final allocation of generations/consumptions, $u^\star$, exists as a limit, and is optimal for (\ref{deterministic-optimal-control},\ref{systemconstraint}).  
If $\lim_ku_i^k$ exists, averaging is not needed since it is equal to $u^\star$.
\hfill $\square$.
\end{theorem}
%
%

In this deterministic context, the whole problem can be solved at time 0, 
with actions $u^{\star} := (u^{\star}(0), u^{\star}(1), \ldots , u^{\star}(T-1))$ implemented open loop.

\section{Iterative Bidding Schemes for Stochastic Systems}\label{sec-ibs}
\begin{figure}
\centering
\includegraphics[width=0.8\linewidth]{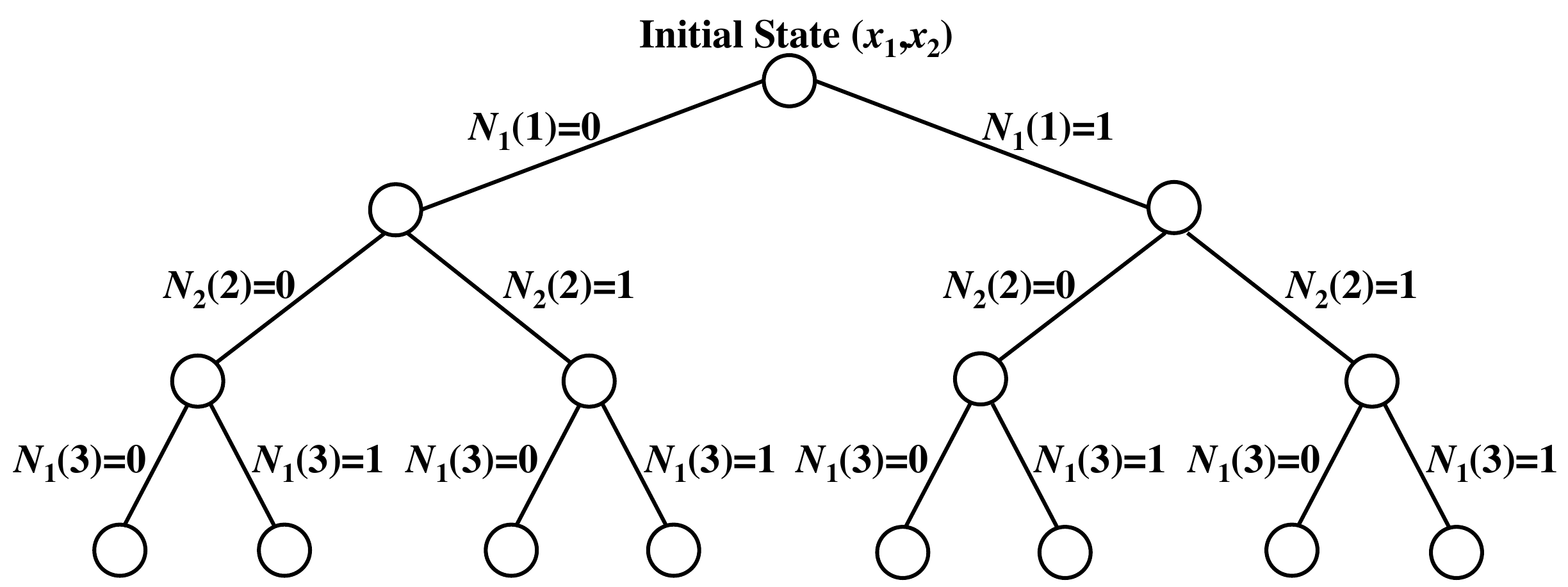}
\caption{A tree visualization of uncertainty for a two agent system evolving over three bid times, where the uncertainty values are binary, either 0 or 1.}
\label{extenform}
\end{figure}

We now turn to the stochastic case.
The goal is to solve the ISO Problem~(\ref{p0}) through Iterative Bidding Schemes (IBS), as in Walrasian tatonnement~\cite{arrow}. We explain what transpires in such an IBS for the simpler common uncertainty context $N(t) \equiv N_c(t)$.
A tree visualization of the system randomness, as in Fig.~\ref{extenform}, is helpful. 
Suppose that $N(t)$ assumes only finitely many values. We can then construct an uncertainty tree of depth
$T$, in which the root node corresponds to the initial system state, and
the sequence of transpired noises $\{ N(0), N(1), \ldots , N(s-1) \}$ corresponds to 
some node at the level $s$. 

Since all agents know the law of $\{ N_c(t) \}$, i.e., the probability measure induced on the sample paths of the noise stochastic process $\{ N_c(t) \}$, the agents know the topology of the tree, and the transition probabilities along edges.
However, the agents do not know the system dynamics of other agents, their utility functions, or states or actions. 
The ISO need not know the law of $N$.
We will suppose that the ISO does know the topology of the tree and the labels of the nodes.

Let $\mathcal{F}_{i,t} := \sigma(X_i(0), N_i(0), N_i(1), \ldots , N_i(t-1), N_c(0), N_c(1), \ldots , N_c(t-1))$ denote the sigma-algebra 
generated by agent $i$'s observations up to time $t$. (Incorporating $N_i (\cdot )$ is unnecessary since
private uncertainties are absent here, but will be useful subsequently in the general case of private observations).

The IBS scheme will intertwine two processes, a \textit{Bid Update Process} and a \textit{Price Update Process}.
As in Section \ref{sec-static}, information revealed during the bidding process is important to determining the final allocation.
Additionally, repeating the process at each time instant is important in the stochastic dynamic case in adapting to how agents are evolving over time as uncertain events happen. \\
{\bf{Bid Update Stochastic Process $\mathbf{B_s = (U_{i,s}(s), U_{i,s}(s+1),\ldots, U_{i,s}(T-1))}$}}: The bid update stochastic process $B_s$ 
\emph{at the particular time} $s$ of an agent $i$ specifies how much electricity that agent intends to purchase (negative if supplying) in every time period from that time $s$ till the final time $T-1$ in response to future events. 
As above, for illustratory purposes, we assume that $N(\cdot)$ is observed
causally by all agents. Then, this bid function of agent $i$ is a function which specifies to the ISO, at any time $s$, as a function of the past history of observed noise $N(\tau), \tau <s$, how much electricity it will purchase at each instant in the future under different future uncertainties. In Fig.~\ref{extenform},
the bid function of agent $i$ specifies, for each node in the tree, the amount of electricity that it is willing to purchase 
if and when the system passes through that node. 

{\bf{The Price Update Stochastic Process $\mathbf{\lambda_s = (\lambda_s(s), \lambda_s(s+1), \ldots , \lambda_s(T-1))}$}} is a
stochastic process announced by the ISO at time $s$. Assuming that the noise process $N( \cdot )$ is observed causally
by all the agents, it specifies for each time $s \leq t \leq T-1$, as a function of the past history of observed noise $N(\tau), \tau < s$, the price $\lambda_s(t)$ at which electricity will be sold in the market
at time $t$ under different future uncertainties. In the tree of Fig.~\ref{extenform}, it corresponds to a price corresponding to each node of the tree at levels $s$ through $T-1$. 

{\bf{$\mathbf{k}$-th Bid Update at time $\mathbf{t}$:}} Suppose that the ISO has declared a price  process $\lambda^k_s$ at time $s$,
where $k$ is an index that
we will use for iteration. In the \textit{Bid Update}, each agent $i$ changes its bid $B^k_s$ in response to the price function $\lambda^k_s$ by solving the following problem, dubbed Agent $i$'s Problem,
\begin{align}
&\min_{U_i \mbox{\scriptsize{ s.t. }} (\ref{constraints-at-time-t-1},\ref{constraints-at-time-t-2})} \mathbb{E} [ \sum_{t=s}^{T-1} [ c_i(X_i(t),U_i(t),t) + e_i(U_i(t),t) \notag \\
& \quad \quad \quad \quad\quad \quad \quad  + \lambda^k_s(t)^T E_i(t) U_i(t) ] |  \mathcal{F}_{i,s} ]. \label{comp-obs-prob}
\end{align}

{\bf{$\mathbf{(k+1)}$-th Price Update at time $\mathbf{s}$:}} The ISO updates the price process in response to the agents' bids. 
Guided by the ``excess consumption function" $\sum_{i=1}^M E_i(t) U^k_{i,s}(t) - d(t)$,
it raises or lowers prices to satisfy the general linear constraint:
\begin{equation}
\lambda^{k+1}_s(t) = \lambda^{k}_s(t) + \frac{1}{k} [\sum_{i=1}^M E_i(t) U^k_{i,s}(t) - d(t)], \quad s \leq t \leq T-1. \label{priceupdate}
\end{equation}

{\bf{The Final Averaged Allocations.}} At time $s$, after the prices have converged, i.e.,
$\lambda^{\star}_s(t) = \lim_{k \to \infty} \lambda^{k}_s(t)$ for $s \leq t \leq T-1$,
the ISO announces the allocations as the limits 
$U^\star_{i,s}(t) := \lim_{k \to \infty}\bar{U}^k_{i,s}(t)$ for $s \leq t \leq T-1$, of  the following averaged
bids,
\begin{align}
&\bar{U}^k_{i,s}(t) = \frac{\sum_{s=1}^{k-1}s^\theta}{\sum_{s=1}^{k}s^\theta} \bar{U}^{k-1}_{i,s}(t) + 
\frac{k^\theta}{\sum_{s=1}^{k}s^\theta} U^{k}_{i,s}(t), \label{averaged-stock-bids}
\end{align}
with $\bar{U}^0_{i,s}(t)= U^0_{i,s}(t)$. If the unaveraged agent bids converge, then their limit is the same as the above.
As in Section~\ref{sec-static}, under appropriate conditions the limits of the prices and
the averaged bids do exist.


\section{Stochastic Systems with Common Uncertainties: State Prices and Bidding} \label{sec-common}
Now we analyze how the above Iterative Bidding Scheme functions in the case of common uncertainties,
i.e., $N(t) \equiv N_c(t)$, and there are no private noises $N_i$.
Denote the combined state of the system by $X(t) := (X_1(t), X_2(t), \ldots , X_M(t))$, and the combined actions by $U(t) := (U_1(t), U_2(t), \ldots , U_M(t))$.
%

At each time $s$, 
a sequence of iterative \emph{tentative} price announcements by the ISO 
for each node at or below the current node at level $s$, followed by \emph{tentative} bids by all agents 
for such nodes responding optimally
to the price announcement, takes place, until they converge. 
At each iteration, the ISO revises the tentative price announcement to drive the ``excess consumption" at each node towards zero, and agents respond optimally according to their own cost-to-go
function. This iteration of tentative
prices and tentative bids continues till the prices converge. At that point 
the agents consume/generate the
weighted average amount they bid for the particular node occupied at time $s$. 
The system then moves forward to time $s+1$, arriving at a
random node at level $s+1$ according to $N(s)$, and the entire process is repeated. This is in the same fashion as Model Predictive Control.

This process is a dynamic modification of Arrow's~\cite{arrow} approach of treating each ``good" available at a certain time and place as a separate good.  Since agents do not know each other's dynamics or states or actions there is the added 
critical proviso of bidding for future ``time-places" by each agent,
with only
the current price being actually implemented a la Model Predictive Control.


\begin{assumption}[Assumption for stochastic case]\label{finite-compact-slater-stochastic}
(i) There is an optimal solution of (\ref{p0}) with finite cost. \\
(ii) $ \sum_{t=0}^{T-1} c_i(X_i(t),U_i(t),t)$, $\sum_{t=0}^{T-1} e_i(U_i(t),t)$,  
and $h_i(N_i^t,N_c^t,t,U_i(t))$
are convex in $U_i^{T-1}$ for each noise sequence
$N^{T-1}$,
with $X_i(t)$ a function of $U_i^{t-1}$ and
$N^t$.
\\
(iii) For each fixed noise sequence $N_c^t,N_i^t$, there exists a feasible $\bar{u}$ 
satisfying $h_i(,t,\bar{u}_i(t)) < 0$ in
$\mbox{RelInt(Dom}(c_i)) \cap \mbox{RelInt(Dom}(e_i))$ for
$\mbox{ for }1 \leq i \leq M, 0 \leq t \leq T-1$.
For simplicity of exposition we suppose that the noise processes $N_c(t), N_i(t)$ assume only finitely many values, allowing them to be represented by a tree as in Fig.~\ref{extenform}. 
\end{assumption}


\begin{theorem}
In the above common uncertainty case, the price-bid solution, with price updates~(\ref{priceupdate}), and bid updates determined as the optimal solution of (\ref{comp-obs-prob}), and allocations at each $t$ given by $\lim_{k \to \infty} \bar{U}^k$
where $\bar{U}^k$ is obtained as the averaged version of $U^k$ as in (\ref{ergodic_method}), achieves the
maximum social welfare when the cost functions satisfy Assumption~\ref{finite-compact-slater-stochastic}.
\end{theorem}
\noindent \emph{Proof:}
Let us suppose that $x(0)$ is fixed, without loss of generality. Let $p_v$ denote the probability of node $v$ in the uncertainty tree. The depth of the node in the tree indicates time. Now, a Markov policy~\cite{kumar} maps the system state and time to actions, thereby specifying an action $U(v) := (U_0(v), U_1(v),\ldots, U_M(v))$ satisfying $\sum_{i=1}^M E_i(v) U_i(v) = d(v)$ for every node $v$ in the tree. This is easily seen to be true by noting that each node in the uncertainty tree also indicates the state of the system at that time. Now let us consider a more general ``\emph{tree policy}" that specifies a $U(v) := (U_0(v), U_1(v),\ldots, U_M(v))$ satisfying $\sum_{i=1}^M E_i(v) U_i(v) = d(v)$ for every node $v$ in the tree. The class of tree policies is more general than the class of Markov policies, since two nodes in the tree at the same depth may correspond to the same state $X(t)$
arrived at through differing noise realizations, but a tree policy can choose different actions for them. Since the class of Markov policies contains an optimal policy, it follows that the class of tree policies also contains one. 

For every tree policy, for every node $v$, there is a unique sequence of actions $U^v := \left\{U(0),\ldots,U(t)\right\}$ taken in the preceding $t$ steps, where $t$ denotes depth of node $v$. The state $X(t)$ corresponding to $v$ is thereby determined by $(v,U^v)$. The problem~(\ref{p0}) can be written equivalently as the following optimization problem,
\begin{align}
& \mbox{Min} \sum_{i=1}^{M}\sum_v p_v [c_i \left(v,U^v\right)+ e_i \left(U^v\right) ] \notag \\
& F_i(v)U_i(v) \leq g_i(v)
+ \sum_{\{v': v' \mbox{ \scriptsize{a predecessor of} } v \}} C_{i}(v',v)U_i(v'), \label{stochastic-nodal-constraints-1} \\
&h_i(U_i(v),v) \leq 0, \label{stochastic-nodal-constraints-2} \\
& \mbox{ such that } \sum_{i=1}^M E_i(v) U_i(v) = d(v),\forall v. \notag
\end{align}
Under Assumption~\ref{finite-compact-slater-stochastic}, the convex programming problem has no duality gap. Let $\lambda(v)$ be the Lagrange multiplier for the constraint $\sum_{i=1}^M E_i(v) U_i(v) = d(v)$, and define the vector $\lambda : = \{\lambda(v) \}$. We obtain,
\begin{align*}  
 \mathcal{L}\left(U, \lambda \right)  := \sum_{i=1}^{M}\sum_v p_v & [ \sum_v c_i(v,u^v) + e_i(u^v) + \\
& \lambda (v)^T E_i(v) U_i(v) ] - \sum_v p_v \lambda (v)^T d(v).
\end{align*}
The process $\lambda(v)$ will be called the ``price process".
Each agent submits a bid for each possible future realization $v$ of the noise process, while the ISO specifies a price at each $v$. The proof parallels the deterministic one.
\hfill $\square$


\section{Stochastic Systems with Private Uncertainties} \label{sec-private}
Now we address the general case where agents have private uncertainties in addition to common uncertainty, i.e., $N(t) = (N_c(t), N_1(t), N_2(t), \ldots , N_M(t))$, where $N_c$ is a common uncertainty that is observed by all, but each $N_i$ is only observed by agent $i$. 
We will suppose that all the agents know the law of $\{ N(t) : 0 \leq t \leq T-1 \}$, but that the ISO knows only the labels of the noise. The agents do not know the dynamics or the cost/utility functions or states of the other agents.

The same Price-Bid iteration can be used, as detailed in Algorithm 
\ref{a2}, and the same result carries over.
\begin{theorem}
For the above system featuring private uncertainties as well as common uncertainties, the bidding process with ISO updating prices according to (\ref{priceupdate}), each agent $i$ updating its consumption/generation bid  according to the optimal solution of (\ref{comp-obs-prob}), and allocations determined at each $t$ by the averaging as in (\ref{ergodic_method}),
achieves the optimal social welfare under Assumption~\ref{finite-compact-slater-stochastic}.
\end{theorem}
\noindent \emph{Proof:}
At time $0$, there are no private noises $N_i(-1)$, and so the above proof holds at time $0$. Noting this, it follows that the result also holds at each time $s \geq 1$ since the bid-price iteration is repeated at each such time, and we can simply regard $s$ as the new ``initial" time. 
\hfill $\square$

\begin{algorithm}
\caption{: Stochastic Dynamic Agents}
\label{a2}
\begin{algorithmic}
\STATE \textbf{Assumption:}  The law of the combined noise process $\mathcal{L}(N)$ is common knowledge of all agents and labels are known to the ISO. 
 \FOR{ bidding times $s=0$ to $T-1$}
\STATE $k=0$
\STATE \REPEAT
\STATE Each agent $i$ solves the problem
\begin{align*}
&\mbox{Min} \quad \mathbb{E}\sum_{t = s}^{T-1} [ c_i(X_i(t),U_i(t),t) + e_i(U_i(t),t) \\
& \quad \quad \quad \quad \quad \quad \quad \quad + \lambda^{k}(t)^T E_i(t) U_i(t) ] ,
\end{align*}
with initial condition $X_i(s)$ to obtain the optimal $\{U^k_{i,s}(t),s\leq t\leq T-1\}$
subject to (\ref{stochastic-nodal-constraints-1},\ref{stochastic-nodal-constraints-2}), and submits it to ISO.
\STATE The ISO declares new prices for $s \leq t \leq T-1$, i.e.,
\begin{align*}
& \lambda^{k+1}(t) = \lambda^k(t) + \frac{1}{k} [\sum_{i=1}^M E_i(t) U^k_i(t) - d(t)].
\end{align*}
\STATE $k\to k+1$
\UNTIL{$\lambda^k(t)$ converges a.s. to $\lambda^{\star}(t) \mbox{ for } s \leq t \leq T-1$.}
\STATE ISO computes $\bar{U}^k_{i,s}(s)$ as in (\ref{averaged-stock-bids}), and implements $U^\star_{i,s}(s) := \lim_{k \to \infty}\bar{U}^k_{i,s}(s)$.
\ENDFOR
\end{algorithmic}
\end{algorithm}


The assumption that the law of $\mathcal{L}(N(t))$ is common knowledge can potentially be relaxed by utilizing Stochastic Approximation~\cite{borkarbook,kushner,robbins}, so that agents can ``learn" them as the system evolves.

The major drawback of this algorithm is that it is exponentially
complex in $T$ due to the number of states in the tree, even if each $N(t)$
is binary. In the next section we show that we can dramatically simplify the bidding process and solution in the LQG context.

\begin{figure}[h]\centering\resizebox{7cm}{10.5cm}{ 
\begin{tikzpicture}[node distance=2cm]
\tikzset{trapezium stretches=true}
\node (start) [startstop] {Start with $s=0$};
\node (in1) [io, below of=start] {ISO declares prices $\lambda^0 (t)$ for times $s \leq t \leq T-1$, sets $k=0$ };
\node (pro1) [process, below of=in1] {Each Agent $i$ optimizes consumption/generation for $s \leq t \leq T-1$
for the \emph{deterministic} model $x_i(t+1) = A_ix_i(t) + B_iu_i^k(t)$, $x_i(s)=X_i(s)$,
with cost $\sum_{t=s}^{T-1} [ x_i^\intercal (t)Q_i x_i(t) + u_i(^kt)^\intercal R_i u_i^k(t) + \lambda^k(t)u_i^k(t)]$ }; 
\node (in2) [io, below of=pro1] {Agents submit bids $u_i^k(t)$ for $s \leq t \leq T-1$}; 
\node (dec1) [decision, below of=in2,yshift=-1cm] {Bids converged?};
\node (pro2) [upd, right of=dec1,xshift=1.8cm] {Update Prices $\lambda^{k+1}(t) = \lambda^k(t) + \alpha^k \sum_i u_i^k(t)$ for $s \leq t \leq T-1$};
\node (pro3) [process, below of=dec1,yshift = -1.2cm] {Implement the first entry of the converged bids as $U_i(s)$. Agents update states $X_i(s+1)=A_iX_i(s) +B_i U_i(s) + N_i(s)$. Increment time $s$ by one};
\node (dec2) [decision, below of=pro3,yshift = -1cm] {Is $s=T$?};
\node (stop) [startstop, below of=dec2] {Stop};
\node (dummy) [circle,radius = .0pt,inner sep=0pt,right of = dec2,xshift = 3cm]{};
\draw [arrow] (start)--(in1);
\draw [arrow] (in1)--(pro1);
\draw [arrow] (pro1)--(in2);
\draw [arrow] (in2)--(dec1);
\draw [arrow] (dec1)-- node[anchor=east] {yes} (pro3);
\draw [arrow] (pro3) -- (dec2);
\draw [arrow] (dec2) --node[anchor=east] {yes}  (stop);
\draw [arrow] (pro2) |- (pro1);
\draw [arrow] (dec1) -- node[anchor=south] {no} (pro2);
\draw [thick] (dec2) -- node[anchor=south]{no} (dummy);
\draw [arrow] (dummy) |- (in1);
\end{tikzpicture}
}
\caption{Scheme for ISO Problem with LQG Agents.}
\label{flo}
\end{figure}
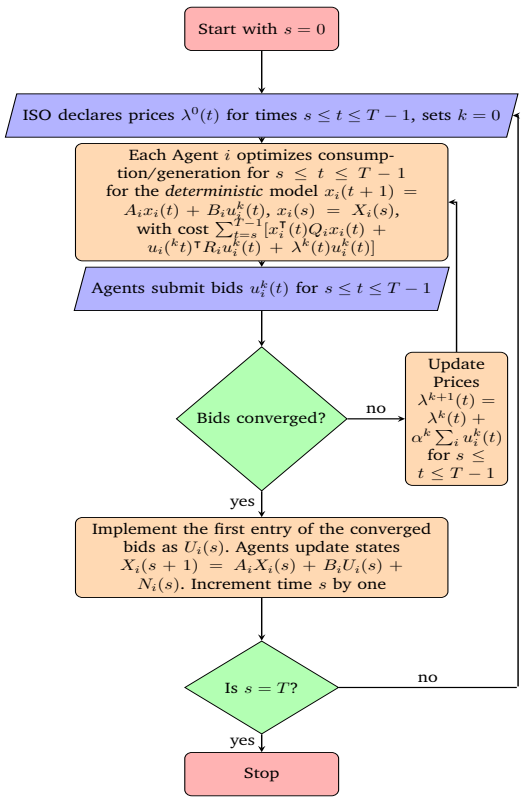

\section{The ISO Problem for LQG Agents} \label{sec-lqg}
We will now show that in the LQG case one can meet both the stringent privacy and lack of knowledge constraints of other agents, and yet avoid the complexity of the solution in the general case where stochastic process bids need to be made for all future uncertainties. Only iterations between \emph{vectors of future prices} announced by the ISO, and \emph{vectors of future consumption/generation bids} by the agents are needed,
similar to deterministic dynamic systems. Moreover, agents need not know the laws of the private uncertainties of other agents
or anything at all about each other. In fact they do not even need to know of the existence of the others. Yet, optimal coordination can be achieved by the ISO, and that too in a tractable manner where agents bid at each time. This appears practically feasible with bid periods separated by minutes.

The only difference between the deterministic 
dynamic case treated in Section \ref{sec-deterministic} and the LQG case is that while in the former the bid-price iteration only needs to be carried out at time $0$, in the LQG case it needs to be carried out at each time $s$. This is similar to Model Predictive Control, where we only implement the first step of the prices and consumptions/generations at each time $s$. 

The $M$ agents have linear dynamics affected by Gaussian noise and have quadratic costs. 
Externalities constituting a positive semidefinite quadratic plus a linear term in $U_i$ could be included, but are omitted below for simplicity.
Initial conditions and noises are Gaussian:
$X_i(0) \sim N(0, \Sigma_{i,0})$ and $N_i(t) \sim N(0, P_{i,t})$, and independent of all others. The cost functions of agents, are quadratic, with
$Q_i \geq 0$ and
$R_i >0$. The ISO Problem is:
\begin{align}
& \mbox{Min} \ \mathbb{E} \sum_{t=0}^{T-1}\sum_{i=1}^{M} [ X^\intercal_i(t)Q_i X_i(t) + U^{\intercal}_i(t) R_i U_i(t) ] 
\label{lqgisocost} \\
&X_i(t+1) = A_iX_i(t) +B_i U_i(t) + N_i(t),
\mbox{ and } (\ref{constraint-on-balancing}) \label{lqgisosystem}.
\end{align}
The case of time-varying systems is entirely analogous.


Agents have no knowledge of each other's presence.
Agent $i$ does not know the value of $M$, the number of agents,
the matrices $\left\{ A_j,B_j,Q_j,R_j, \Sigma_{j,0}, P_{j, \cdot} ) \right\}_{j\neq i}$
of other agents, the realizations of their state processes $\{ X_j (\cdot ), j\neq i \}$ or noises $\{ N_j( \cdot ), j\neq i \}$.

We propose and prove the convergence
and optimality of an Iterative Bidding Scheme which is much simpler than that of Section~\ref{sec-common} in the following critical aspect. The bid function submitted at time $s$ specifying the quantity of electricity that agent $i$ is willing to purchase at times $\{ s,s+1,\ldots,T-1 \}$ is \emph{not} a function of the outcomes of 
the noise sequence $\{ N(t),t>s \}$. It is simply a vector $(u_i(s),u_i(s+1),\ldots,u_i(T-1))$ comprised of $T-s+1$ entries. The same is also true for prices.
The ISO just specifies a vector $(\lambda(s), \lambda(s+1), \ldots , \lambda(T-1))$ of $T-s$ entries. Both are \emph{not} specified as functions of future uncertainties.

This makes it different from
Arrow~\cite{arrow}.
The complexity of specifying prices or bids for all future events is avoided.
Removal of future event-based bidding and prices leads to a drastic reduction in 
the complexity of the iterative scheme that arises
even if all uncertainties were finite valued, let alone real valued as here. 

Another simplifying feature is that the ISO need not average the bids as in the (\ref{ergodic_method}).
The bids of agents converge at each time instant without averaging.

Note that even though the bid function at each time $s$ is not future event-based, it is determined afresh at each time.
At each time $s$, the following iteration takes place: Each agent bids a vector of future generations/consumptions in response to prices announced by the ISO for future power,
and the ISO updates the prices in return, until convergence.
Hence, the converged prices and consumptions/generations do depend on the system states of the $M$ agents, and are therefore stochastic. 

The key to showing the existence of such a simple bidding scheme lies in utilizing the certainty equivalence property of LQG systems~\cite{kumar}.
\begin{algorithm}
\caption{: ISO Problem with LQG Agents}
\label{a3}
\begin{algorithmic}
 \FOR{bidding times $s=0$ to $T-1$}
\STATE $k=0$
\STATE  Initialize $\{ \lambda^k_s (t) : s \leq t \leq T-1 \}$ arbitrarily.
 \STATE \REPEAT
 \STATE Each agent $i$ solves the problem (\ref{agentibidlqrpricecosr})
 for a \emph{deterministic} system (\ref{agentibidlqrsystem})
with initial condition $x_i(s) := X_i(s)$, where $X_i(s)$ is the state of the $i$-th agent at time $s$,
and submits the optimal values, denoted $u_{i,s}^k(t)$, for $s \leq t \leq T-1$ to the ISO.
 \STATE ISO updates the prices according to (\ref{isolqrprice2},\ref{stepsize}). \\
Increment $k$ by $1$.
\UNTIL{$u^k_{i,s}(t)$ converges to $u^{\star}_{i,s}(t)$},
\STATE Implement $( u^\star_{1,s}(s), u^\star_{2,s}(s), \ldots , u^\star_{M,s}(s) )$
\ENDFOR
 \end{algorithmic}
 \end{algorithm}
 
The iterative bidding scheme is shown in Fig.~\ref{flo} and Algorithm \ref{a3}.
For simplicity, consider only balance of energy.
At time $s$, in response to the $(k-1)$-th iterate  announced by the ISO of the price sequence $(\lambda^k_s(s),\lambda^k_s(s+1),\ldots,\lambda^k_s(T))$, agent $i$ announces the \emph{optimal open loop sequence} $(u^k_{i,s}(s), u^k_{i,s}(s+1),\ldots,u^k_{i,s}(T-1))$ for the following \emph{deterministic} Linear Quadratic Regulator (LQR) problem:
\begin{align}
&\min \sum_{t=s}^{T-1} [ x_i^\intercal (t)Q_i x_i(t) + u_i(t)^\intercal R_i u_i(t) + \lambda^k_s(t)u_i(t) ]
\label{agentibidlqrpricecosr} \\
&\mbox{s.t. } x_i(t+1) = A_ix_i(t) + B_iu_i(t), \quad s \leq t \leq T-1, \label{agentibidlqrsystem} \\
&\mbox{with initial condition } x_i(s) := X_i(s). \label{agentibidlqrinitial}
\end{align}

The price adjustment now is just for a vector of real numbers at each time $s$:
\begin{align}
& \lambda^{k+1}_s(t) = \lambda^k_s(t) + \alpha^k \sum_{i=1}^M u^k_{i,s}(t), 
s \leq t \leq T-1. \label{isolqrprice2} \\
&\mbox{ where } \alpha^k > 0, \lim_{k} \alpha^k =0, \sum_{k=0}^{\infty} \alpha^k = + \infty. \label{stepsize}
\end{align}

At time $s$, the iterations in $k$ are continued till the
price iterations 
$(\lambda^{k}_s(s),\lambda^{k}_s(s+1),\ldots,\lambda^{k}_s(T-1))$
converge to $(\lambda^{\star}_s(s),\lambda^{\star}_s(s+1),\ldots,\lambda^{\star}_s(T-1))$.
Denote the corresponding limit of the input sequence of agent $i$ by
$(u^{\star}_{i,s}(s),u^{\star}_{i,s} (s+1),\ldots,u^{\star}_{i,s} (T-1))$.
The price at time $s$ is then set to $\lambda^{\star}_s(s)$ and each
agent $i$ applies the input $u^{\star}_{i,s}(s)$. This is repeated at time $s+1$.
\begin{theorem}
The bid-price iteration scheme (\ref{agentibidlqrpricecosr},\ref{agentibidlqrsystem},\ref{agentibidlqrinitial},\ref{isolqrprice2},\ref{stepsize}) achieves the optimal social welfare for the LQG ISO Problem
(\ref{lqgisocost},\ref{lqgisosystem},\ref{constraint-on-balancing}).

\noindent \emph{Proof:}
Let 
\begin{align*}
& x:=(x_1,x_2,\ldots,x_M), u:=(u_1,u_2,\ldots,u_M),\\
& A := diag(A_1, A_2,\ldots,A_M), B:=diag(B_1,B_2,\ldots,B_M),\\ 
& Q=diag(Q_1,Q_2,\ldots,Q_M), R=diag(R_1,R_2,\ldots,R_M),
\end{align*}
 and consider the following \emph{deterministic constrained} LQR problem, with no noise, and featuring energy balance:
\begin{align}
& \min \sum_{t=0}^{T} [x^\intercal (t)Qx(t) + u^{\intercal}(t) R u(t) ]  \label{fullcost:eq} \\
&\mbox{with } x(t+1) = Ax(t) +B u(t), \mbox{ and } (\ref{constraint-on-balancing}). \label{s:eq}
\end{align}

Since the state is affine in $u$, after substituting for the states,
we have a positive definite quadratic programming
problem with equality constraints.
The Karush-Kuhn-Tucker matrix is nonsingular (Section 10.1 of  \cite{boydbook}) since $R_i > 0$,
and so there are unique $u^\star,\lambda^\star$ optimal for the primal and dual, respectively
The Dual function is a differentiable concave quadratic function, and the subgradient method is actually a gradient
method that converges under non-summability of step-sizes, without even
requiring square summability (Section 2.5 of \cite{ber87}).
The bids $u_i^k$ are affine functions of the prices $\lambda^k$. Since prices
satisfy balancing, so does their limit.
Hence this deterministic problem can be solved by the Bid-Price iteration (\ref{adjoinedcost},\ref{PriceIterate})
between the agents and the ISO, as shown for the deterministic problem in Section \ref{sec-static},
to obtain the optimal inputs $u(t)$ for $0 \leq t \leq T-1$. 

However, at the particular time $s=0$ with $x_i(0)=X_i(0)$, the Bid-Price Iteration
(\ref{agentibidlqrpricecosr},\ref{agentibidlqrsystem},\ref{agentibidlqrinitial}) and
(\ref{isolqrprice2}) corresponds exactly to the same Bid-Price Iteration 
(\ref{adjoinedcost},\ref{systemconstraint}) and (\ref{PriceIterate}) as in Section 
\ref{sec-static}. Hence the end result of Algorithm \ref{a3} at time $s=0$ is the optimal action
for (\ref{fullcost:eq},\ref{s:eq}),
\begin{equation}
u(0) = (u_1(0),u_2(0),\ldots,u_M(0)). \label{initial action}
\end{equation}

Now note that due to energy balance, no matter how the first $(M-1)$ agents choose their consumptions/generations, agent $M$'s choice is
forced to be
\begin{equation}
u_M(t) = - \sum_{i=0}^{M-1}u_i(t) \mbox{ for all } t,  \label{solvedconstraint}
\end{equation}
due to the energy balance constraint.
Hence one can substitute for $u_M(t)$ and obtain an equivalent standard, i.e., \emph{unconstrained}, deterministic LQR problem featuring only $(M-1)$ inputs $u_{\scriptsize \mbox{reduced}} := (u_1, u_2, \ldots , u_{M-1})$, where there is no energy balance constraint:
\begin{align}
& \min \sum_{t=0}^{T} [x^\intercal (t)Qx(t) + u_{\scriptsize \mbox{reduced}}^{\intercal}(t) R_{\scriptsize \mbox{reduced}} u_{\scriptsize \mbox{reduced}}(t) ] \label{reducedcost:eq} \\
&\mbox{subject to } x(t+1) = Ax(t) +B_{\scriptsize \mbox{reduced}} u(t), \label{reduceds:eq}
\end{align}
the \emph{deterministic reduced unconstrained} LQR problem.

For this problem
(\ref{reducedcost:eq},\ref{reduceds:eq}), 
which is just a standard LQR Problem, the optimal solution is given by linear feedback
$u_{\scriptsize \mbox{reduced}}(0) = \Gamma_{\scriptsize \mbox{reduced}}(0) x(0)$,
where $\Gamma_{\scriptsize \mbox{reduced}}(\cdot)$ is the optimal feedback gain.

Noting that $u_M$ is linear in $u_{\scriptsize \mbox{reduced}}$, we deduce
that for the full system (\ref{fullcost:eq},\ref{s:eq}) with all $M$ agents, the optimal solution for the
deterministic constrained LQR problem
with the energy balance constraint, is
$u(0) = \Gamma(0) x(0)$,
where $\Gamma(\cdot)$ is the optimal feedback gain obtained from $\Gamma_{\scriptsize \mbox{reduced}}$ through (\ref{solvedconstraint}).

Now consider the corresponding \emph{reduced unconstrained stochastic} LQG problem where there is white Gaussian noise in the state equations~(\ref{s:eq}):
\begin{align}
& \min E \sum_{t=0}^{T} [X^\intercal (t)QX(t) + U_{\scriptsize \mbox{reduced}}^{\intercal}(t) R_{\scriptsize \mbox{reduced}} U_{\scriptsize \mbox{reduced}}(t) ]  \label{stochreducedcost:eq} \\
&\mbox{with } X(t+1) = AX(t) +B_{\scriptsize \mbox{reduced}} U_{\scriptsize \mbox{reduced}}(t) + N(t).\label{stochreduceds:eq}
\end{align}
By Certainty Equivalence~\cite{kumar}, the same linear feedback \emph{gain} as in the deterministic reduced LQR problem is also optimal. In particular, in state $X(0)=x(0)$ at time $0$,
$U(0)=\Gamma(0)x(0)$
continues to be optimal. Thus $u(0)$ given by (\ref{initial action}) is optimal for (\ref{stochreducedcost:eq},\ref{stochreduceds:eq}).

However, reduced unconstrained \emph{stochastic} LQG problem (\ref{stochreducedcost:eq},\ref{stochreduceds:eq}) is equivalent to unreduced constrained LQG problem (\ref{lqgisocost},\ref{lqgisosystem},\ref{constraint-on-balancing}), and so the same $U(0)$ is optimal.

Thus the Bid-Price iteration scheme determines the optimal actions for the agents at time $0$. Our scheme 
(\ref{agentibidlqrpricecosr},\ref{agentibidlqrsystem},\ref{agentibidlqrinitial}) for the LQG problem repeats such a Bid-Price scheme iteration at each time $s=0,1,\ldots,T-1$. Each $X(s)$ can be regarded as an initial state for a subsequent system re-started at time $s$, and the above argument shows that the actions $U(s)$ that it results in for the agents at all times $s$ are also optimal, completing the proof.
\hfill $\square$
\end{theorem}

This result extends to LQG systems where each agent $i$ only has noisy observations $Y_i(t) = D_i X_i(t) + V_i(t)$, where $V_i$ are independent and Gaussian.


\section{Incorporating Additional Linear Constraints: The DC Optimal Power Flow} \label{sec-powerflow}
Besides energy balance,
there are additional constraints of interest. 
An important one is ensuring that the power flows are delivered over the network.
 These constraints are captured by the AC Power Flow Equations~\cite{bergen},
 an
approximation of which 
leads to the so-called DC Power Flow equations that are linear constraints~\cite{bergen}.
The bid-price iterations can be extended to encompass any such additional linear constraints.
The only difference is that there are several prices, one for each constraint,
that each agent needs to incorporate in choosing its actions.

\begin{theorem}
Consider a system consisting of $M$ agents, where each Agent $i$'s system is a Linear
Gaussian System:
\begin{align*}
&X_i(t+1) = A_iX_i(t) +B_i U_i(t) + N_i(t).
\end{align*}
Agent $i$ has a quadratic cost (negative utility):
\begin{align*}
& \min \mathbb{E} \left(\sum_{t=0}^{T-1} [ X^\intercal_i(t)Q_i X_i(t) + U^{\intercal}_i(t) R_i U_i(t) ] \right) .
\end{align*}
There are $N$ linear constraints that need to be satisfied:
\begin{align*}
&\sum_{i=1}^M \gamma_{i,n} U_i(t) =0 \mbox{ for } 1 \leq n \leq N, t=0,1,\ldots,T-1.
\end{align*}

Neither ISO nor agents know the number $M$ or the dynamics/costs/law/states/noises
of other agents.

Consider the following Bid-Multiple Price Iteration.
At each time $s=0,1,\ldots,T-1$, 
at each iterate $k$, in response to prices $\{ \lambda^k_{n,s}(t): s \leq t \leq T-1 \}$, announced by the ISO, agent $i$ solves the deterministic LQR problem:
\begin{align*}
&\min \sum_{t=s}^{T-1} [ x_i^\intercal (t)Q_i x_i(t) + u_i(t)^\intercal R_i u_i(t) + \sum_{n=1}^N\lambda^k_{n,s}(t)u_i(t) ],
\end{align*}
with $x_i(s) = X_i(s)$, determines the optimal $\{ u^k_{s}(t): s \leq t \leq T-1 \}$, and communicates this
sequence to the ISO.
Upon receiving the bids at iterate $k$ from all the agents at time $s$, the ISO updates the $N$
price sequences:
\begin{align*}
& \lambda^{k+1}_{n,s}(t) = \lambda^k_{n,s}(t) + \alpha^k \left(\sum_{i=1}^M \gamma_{i,n} u^k_{i,s}(t)\right),
\end{align*}
for $1\leq n \leq N$ and $s \leq t \leq T-1$, with the step-sizes satisfying (\ref{stepsize}). 
The multiple iterations converge, and let $\{ \lambda^{\star}_{n,s}(t): s \leq t \leq T-1 \}$
denote the limit. Correspondingly let $\{ u^{\star}_{s}(t): s \leq t \leq T-1 \}$ denote the limits of the
bids by the agents.
At each time $s$, agent $i$ applies $U_i(s) = u^{\star}_{s}(s)$.
Then this Bid-Multiple Price Iteration yields the maximum social welfare under the multiple constraints. \\
\noindent \emph{Proof:}
The proof parallel the single constraint case.
\hfill $\square$
\end{theorem}

In the case of the DC Power Flow constraints, this yields the optimal \emph{stochastic dynamic} location marginal prices \cite{bohn1984optimal} that simultaneously take into account all the factors of location,
dynamics and stochasticity.  

\section{Simulation Examples} \label{sec-simu}
In the following, we use the space conditioning example from~\cite{schwepe2} for thermal inertial load agents.
Let $S_1,S_2,S_3$ be sets of conditioning facilities (loads), conventional generators, and renewable suppliers,
respectively, and let $i \in
S_1, j \in S_2, k \in S_3$.
The dynamics of the temperature $X_i(t)$ of
the $i$-th facility is given by (\ref{facility}),
where $X^{O}(t)$ = the outside temperature at time $t$, $\epsilon=\euler^{-\tau/TC}$ = ``factor of inertia", TC = 2.5 hours = time-constant of the system, $\tau$ = time duration between control epochs, which is the same as the inter-bid duration, $\eta =2.5$ = thermal conversion efficiency, and $A=0.14 kW/\degree F$ = overall thermal conductivity.
With $X_i^d(t)$ the desired facility temperature, the cost incurred is a quadratic in the
temperature deviation.
For fossil-fuel generators, the unit-time conventional generation cost curves \cite{wood} for 
supplying energy are quadratic in generation $U_j$. 
We replace hard constraints on ramp-rates $|U_j(t)-U_j(t-1)|$ by a quadratic penalty, with $C_3$ below chosen so that the hard bounds are met, with state given by (\ref{convene}).
For a renewable energy facility $k$, $B_k$ denotes its buffer capacity, $W_k(t)$ stochastic wind/solar energy, and $X_k(t)$ the renewable energy level satisfying (\ref{renlevel}).
Its operating cost is constant.
The resulting ISO Problem~(\ref{p0}) is
\begin{align}
& \min \mathbb{E}\left\{ \sum_{i\in S_1} \sum_{t=0}^{T-1} \left(X_i(t)-X_i^d(t)\right)^2 \right.\notag \\
& \left.+\sum_{j\in S_2}\left( C_{j,1} U_j(t)+C_{j,2} U_j^2(t)+ C_{j,3} \left(U_j(t)-X_j(t)\right)^2	\right) \right\} \notag \\
&\mbox{such that } \sum_{\ell \in S_1 \cup S_2 \cup S_3}^M U_\ell(t) = 0, \mbox{ for }t=1,2,\ldots,T-1, \notag \\
& X_i(t+1) = \epsilon X_i(t) + (1-\epsilon)\left (X_i^{O}(t)+\frac{\eta}{A}U_i(t) \right ), \label{facility} \\
& X_j(t+1) = U_j(t), \label{convene} \\
&   X_k(t+1) = 
\mbox{Min} \{X_k(t)-U_k(t)+W_k(t) ,B_k \} \label{renlevel}. 
\end{align} 

We will compare the performance of the proposed Stochastic Dynamic Optimal Bid-Price Iteration scheme of Sections \ref{sec-private} or \ref{sec-lqg}, called "Optimal" below, with the currently followed Static Dispatch scheme of Section~\ref{sec-static} used in dynamic situations as explained in Section~\ref{sec-intro}, under which the agents perform separate and uncoupled bid-price iterations at each time $t$ to optimize the static cost $C(X(t),U(t))$ incurred at that time $t$. 


\noindent
{\bf{Bidding with LQG Systems:}} A day is divided into twelve $\tau=2$ hour bid-slots, so $\epsilon= 0.4493$. There are only thermal loads,
and wind-farms which have a cost function $\frac{1}{2}X^2(t)$ and with infinite storage capacity $B$. Outside temperatures and 
available wind power are modeled as i.i.d. and normal. (This is only a first step towards modeling the uncertainty,
and other types of distributions can potentially be similarly explored). 
Variance of wind energy is 1 unit for all $t$.
The scenario is described in Table~\ref{table2}.
At the beginning of day, the thermal loads have temperature of $70\degree F$, while wind-farms have 100 units of energy.
The price vector is projected at each update onto a large compact set, and, at termination, the bid vector is projected onto the hyperplane $\sum_i U_i = 0$.

\begin{figure}[h]
\vspace{-0.2in}
\centering
\includegraphics[width=0.75\linewidth]{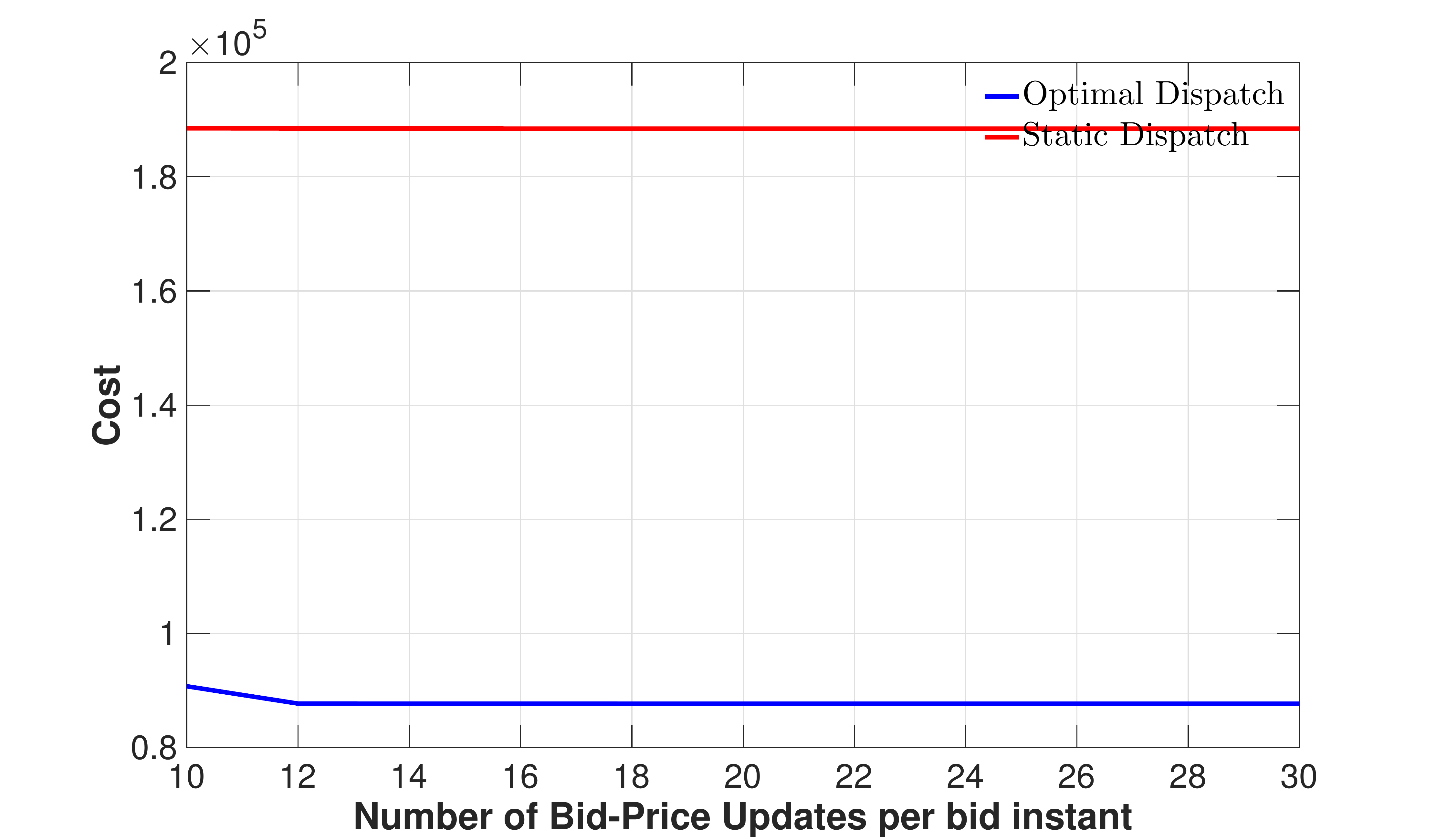}
\vspace{0in}
\caption{Cost, i.e., negative social welfare, vs. number of Bid-Price iterations with five thermal loads and two windfarms.}
\label{fig11}
\end{figure}

\begin{figure}[h]
\vspace{-0.3in}
\centering
\includegraphics[width=0.75\linewidth]{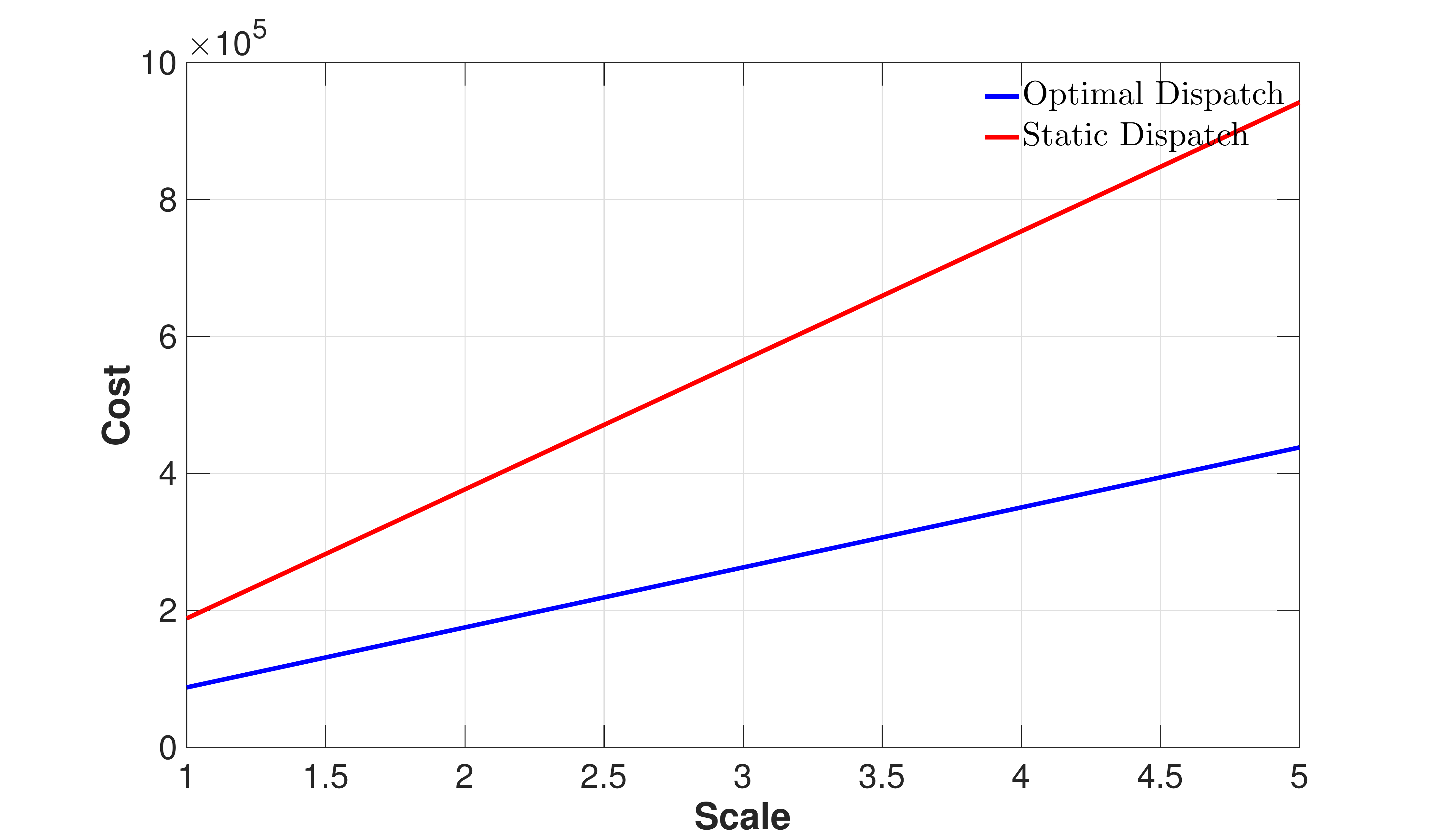}
\vspace{0in}
\caption{Cost as number of users is scaled linearly by $i$, with ratio of thermal loads to windfarms held constant at $5/2$, and with $=15+5 i$ bid-price iterations at each time $t$.}
\label{fig12}
\end{figure}

\begin{figure}[h]
\vspace{-0.1in}
\centering
\includegraphics[width=0.75\linewidth]{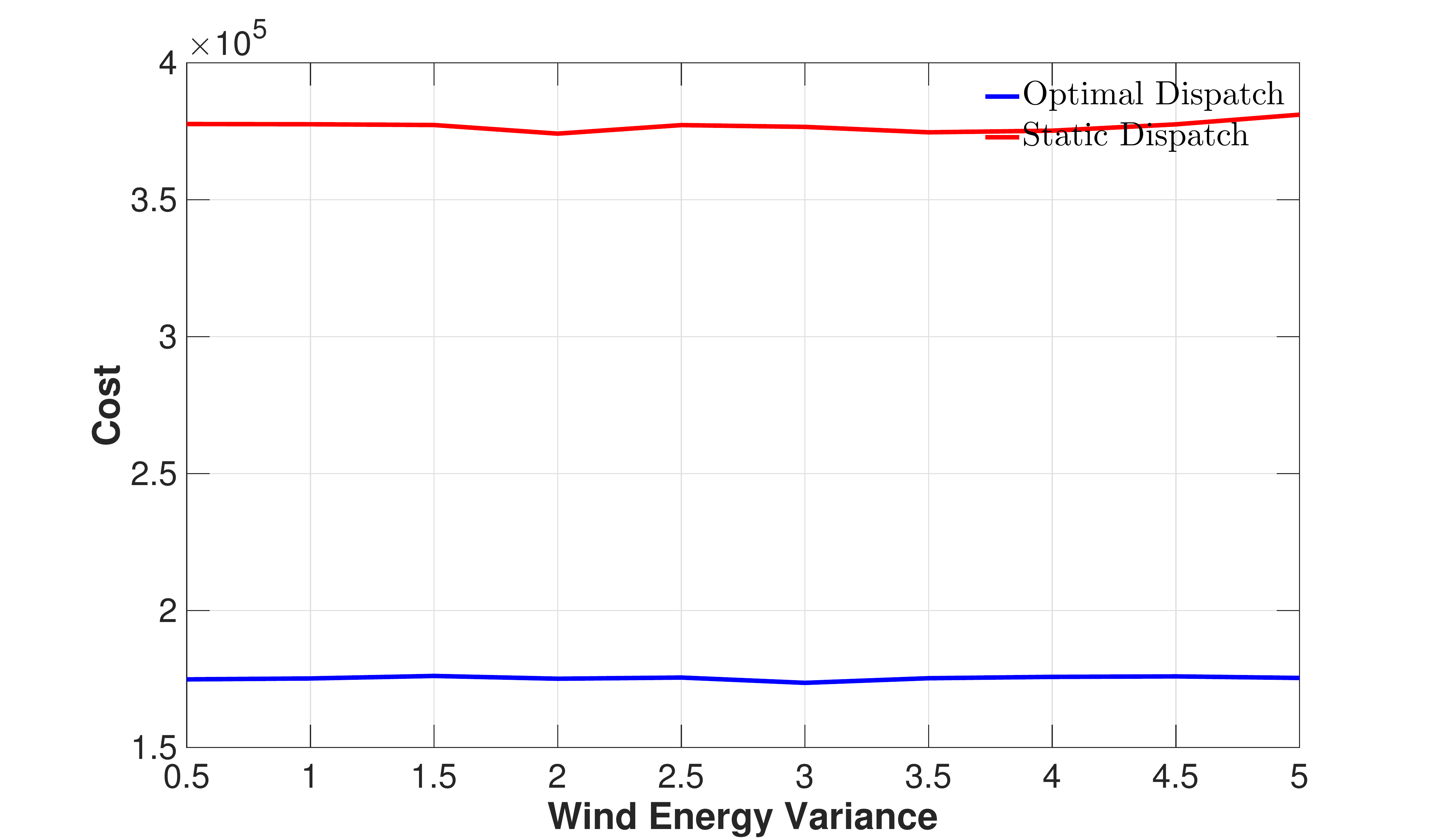}
\vspace{0in}
\caption{Cost vs. wind variance with five thermal loads and two windfarms, with 30 Bid-Price iterations.}
\label{fig13}
\end{figure}

\begin{figure}[h]
\vspace{-0.1in}
\centering
\includegraphics[width=0.75\linewidth]{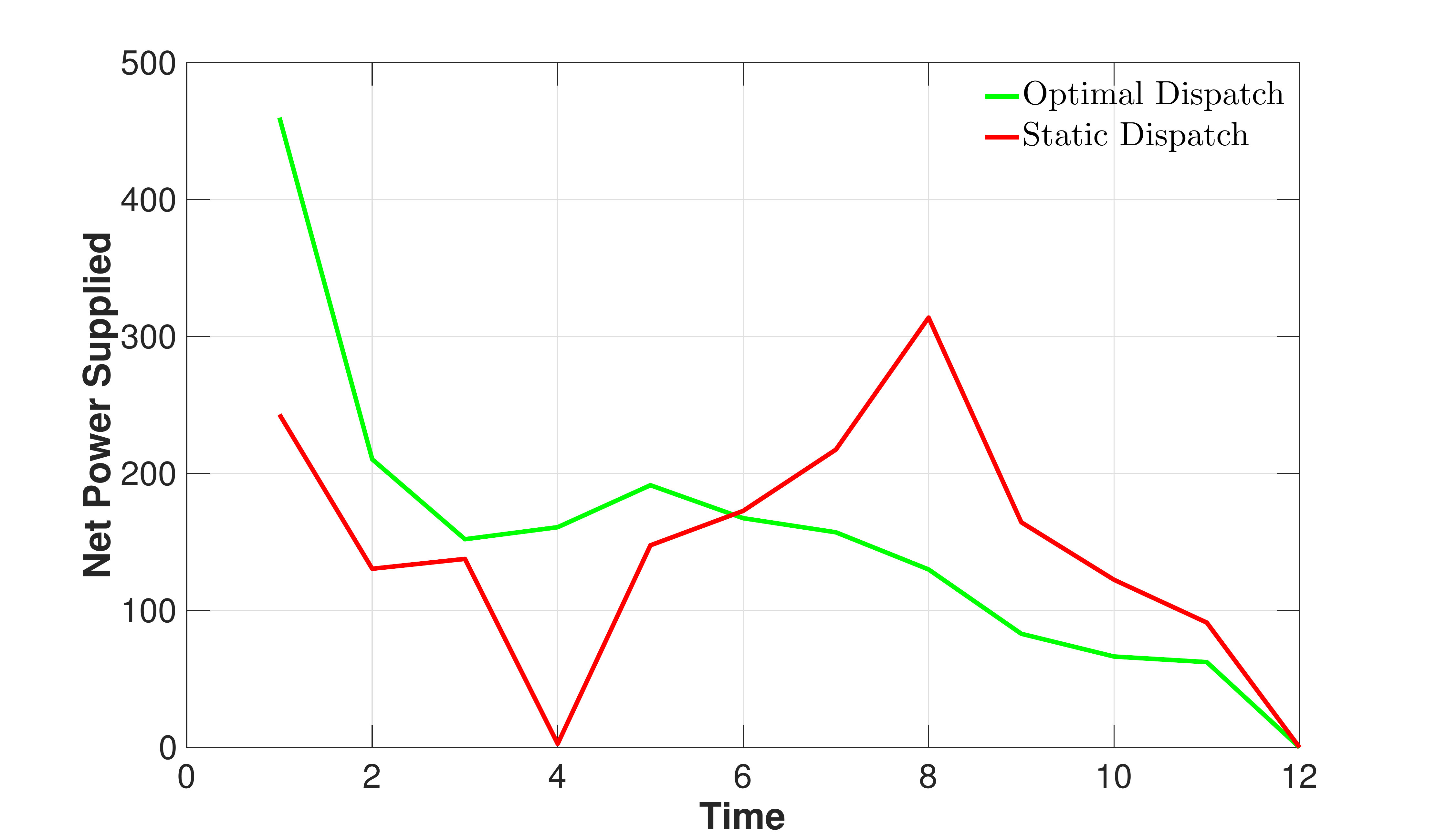}
\vspace{0in}
\caption{Power generation: Optimal scheme ``predicts" the incumbent energy shortage in advance, thereby eliciting smoother generator response. The power production costs for the two schemes are, respectively, $4.37,28.06$ ($\times10^4$), while thermal loads disutility are $13.15,9.62$ ($\times10^4$). Adding these two costs, the net costs are $1.75,3.76$ ($\times10^5$), so that savings achieved by Optimal Scheme is $53.5\%$.}
\label{fig14}
\end{figure}

\begin{figure}[h]
\vspace{-0.1in}
\centering
\includegraphics[width=0.75\linewidth]{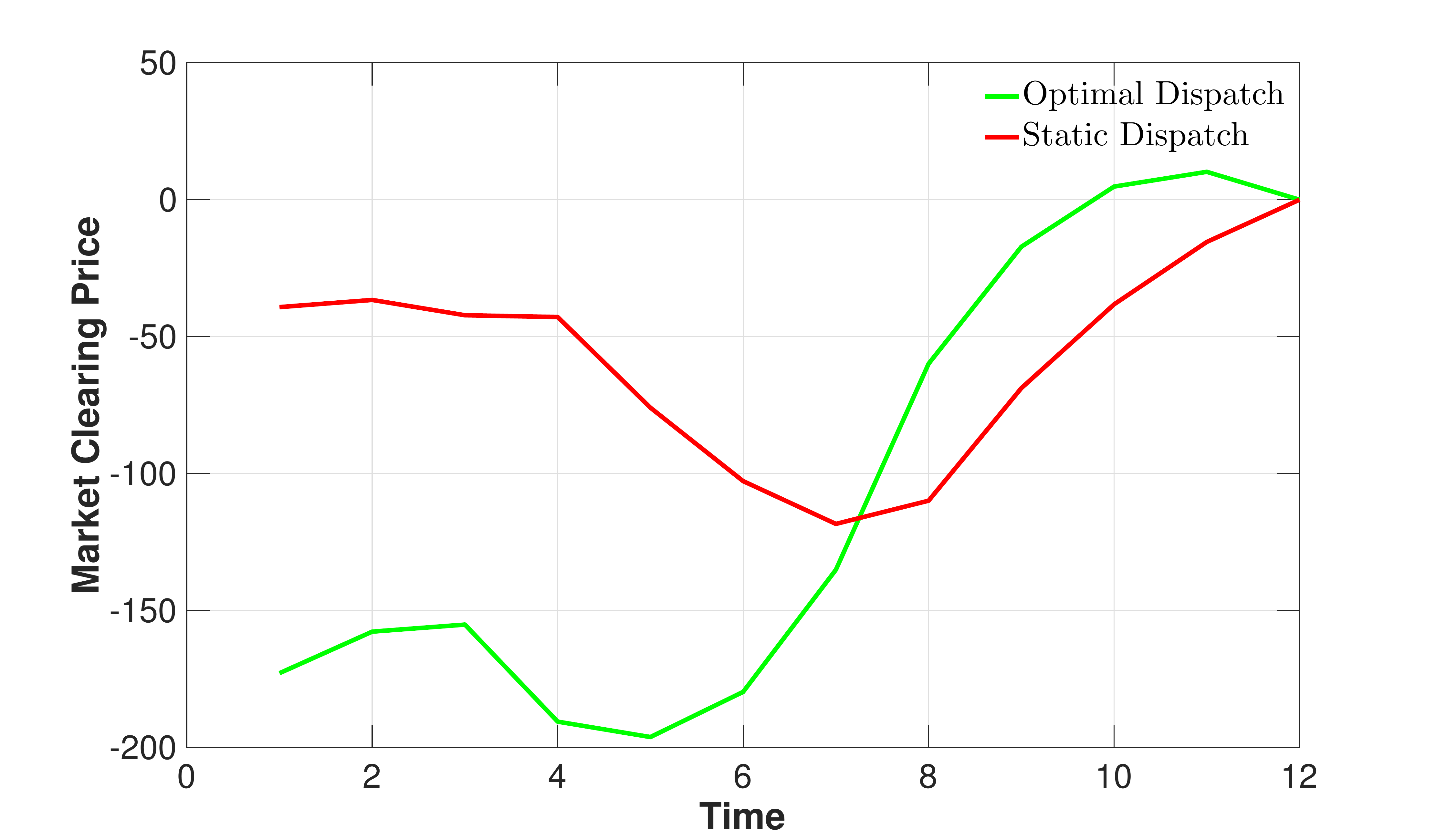}
\vspace{0in}
\caption{Prices: Optimal scheme ``declares" the energy shortage/surplus well in advance, allowing users to react appropriately and eliciting demand response.}
\label{fig15}
\end{figure}

Figs.~\ref{fig11}-\ref{fig13} compare performance of the two schemes as the number of bid-price updates, the number of agents connected to the grid, and variance of wind energy process, are varied. Figs.~\ref{fig14} and \ref{fig15} show how the Optimal scheme is able to attain better social welfare for scale 2.
\begin{table}[h]
\centering
\includegraphics[width=0.75\linewidth]{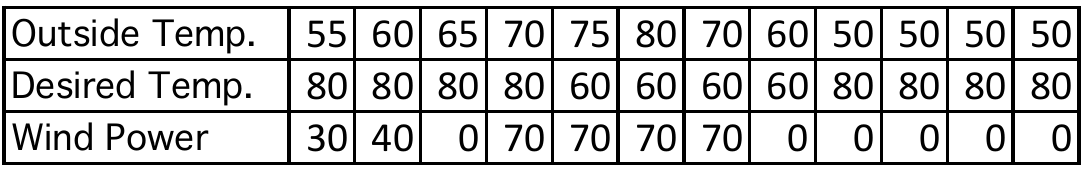}
\caption{Mean outside and desired thermal load temperatures (in $\degree F$), and mean wind power for the 12 periods.}
\label{table2}
\end{table}
\empty
\noindent
{\bf{Bidding in Tree Scenario:}}
The time-horizon is $2$ and time duration between two bids is $5$ hours, roughly coinciding with morning ($7$ am)/$12$ noon, giving $\epsilon = 0.1353$.
Table~\ref{table1} lists stochasticity parameters of wind for two scenarios. 
For fossil plants, $C_1=0.1,C_2=0.01,C_3=0.1$.
Windfarms incur no operational cost.
Bid/price vectors at $t=0$ have three entries, while they are scalar at time $t=1$. 
\begin{table}[h]
\resizebox{8.5cm}{.65cm}{
\begin{tabular}{|l|l|l|l|l|l|l|}
\hline 
   & $\mathbf{T^O(1),T^O(2)}$ & $\mathbf{T^d(1),T^d(2)}$& $\mathbf{W_1,W_2}$ & $\mathbf{P_1,P_2}$ & $\mathbf{S_1/S_2/S_3}$& $\mathbf{B}$ \\\hline
$\mathbf{Case}$ $\mathbf{1}$ & $\mathbf{30,40}$ $\mathbf{(in}$ $\mathbf{\degree F)}$ & $\mathbf{60,80}$& $\mathbf{5,0}$ & $\mathbf{0.5,0.5}$& $\mathbf{7/1/1}$&$\mathbf{30}$   \\\hline
$\mathbf{Case}$ $\mathbf{2}$ &  $\mathbf{40,60}$  & $\mathbf{60,90}$& $\mathbf{10,0}$ & $\mathbf{0.95,0.05}$& $\mathbf{4/1/1}$&$\mathbf{40}$\\\hline
\end{tabular}
}
\caption{
The only stochasticity is wind availability at time $1$, with possible realizations
$W_1,W_2$ with respective probabilities $P_1,P_2$. 
$|S_1|/|S_2|/|S_3|$ are the relative numbers of thermal loads, fossil plants and windmills. }
\label{table1}
\end{table}

\begin{figure}[h]
\vspace{-0.1in}
\centering
\includegraphics[width=0.75\linewidth]{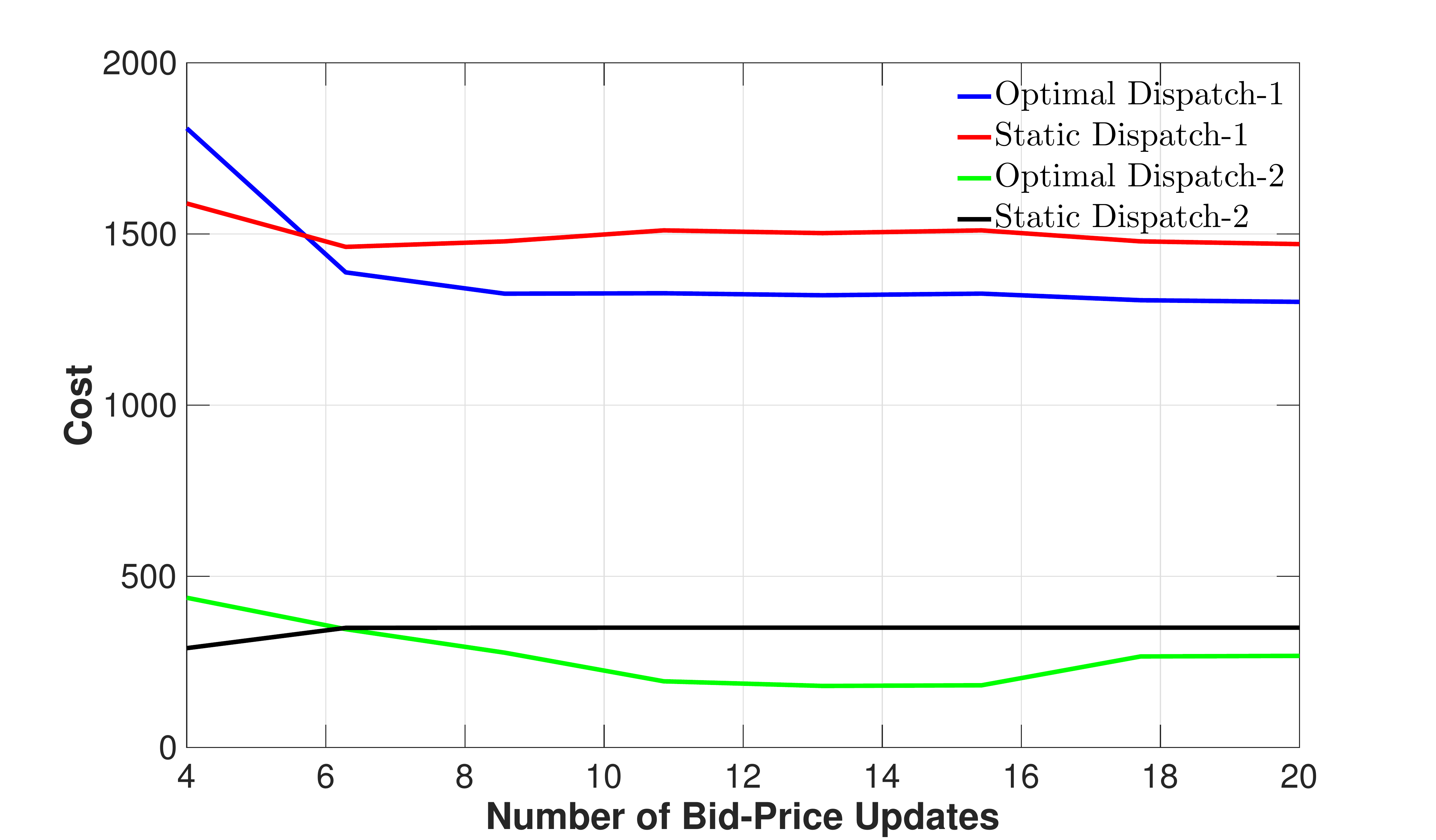}
\vspace{0in}
\caption{Performance as a function of number of Bid-Price updates for the two scenarios in Table~\ref{table1}.}
\label{fig8}
\end{figure}

\begin{figure}[h]
\vspace{-0.1in}
\centering
\includegraphics[width=0.75\linewidth]{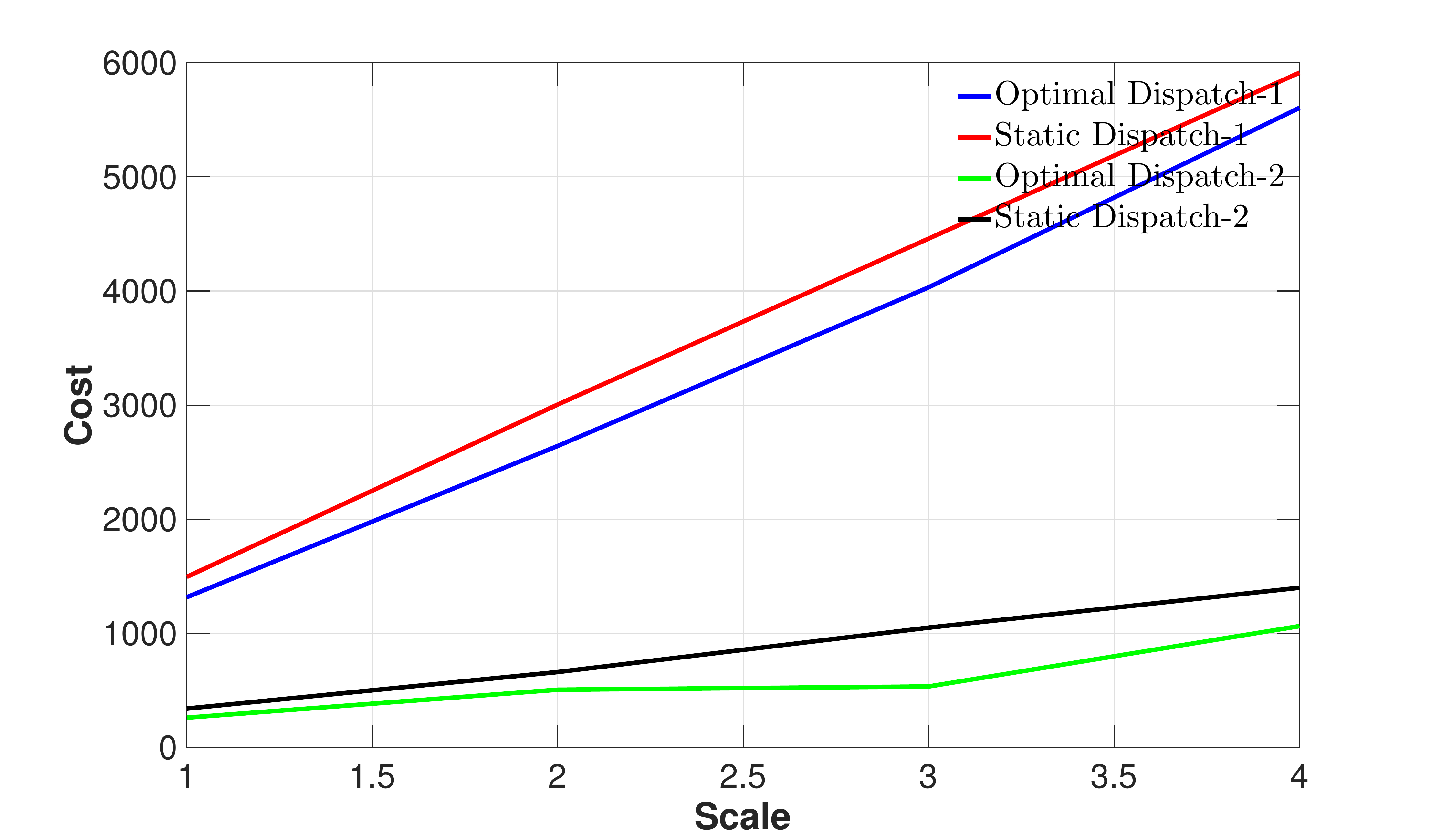}
\vspace{0in}
\caption{Cost as number of agents is increased linearly with scale, in the ratio $S_1/S_2/S_3$ shown in Table~\ref{table1}.}
\label{fig9}
\end{figure}

\begin{figure}[h]
\vspace{-0.1in}
\centering
\includegraphics[width=0.75\linewidth]{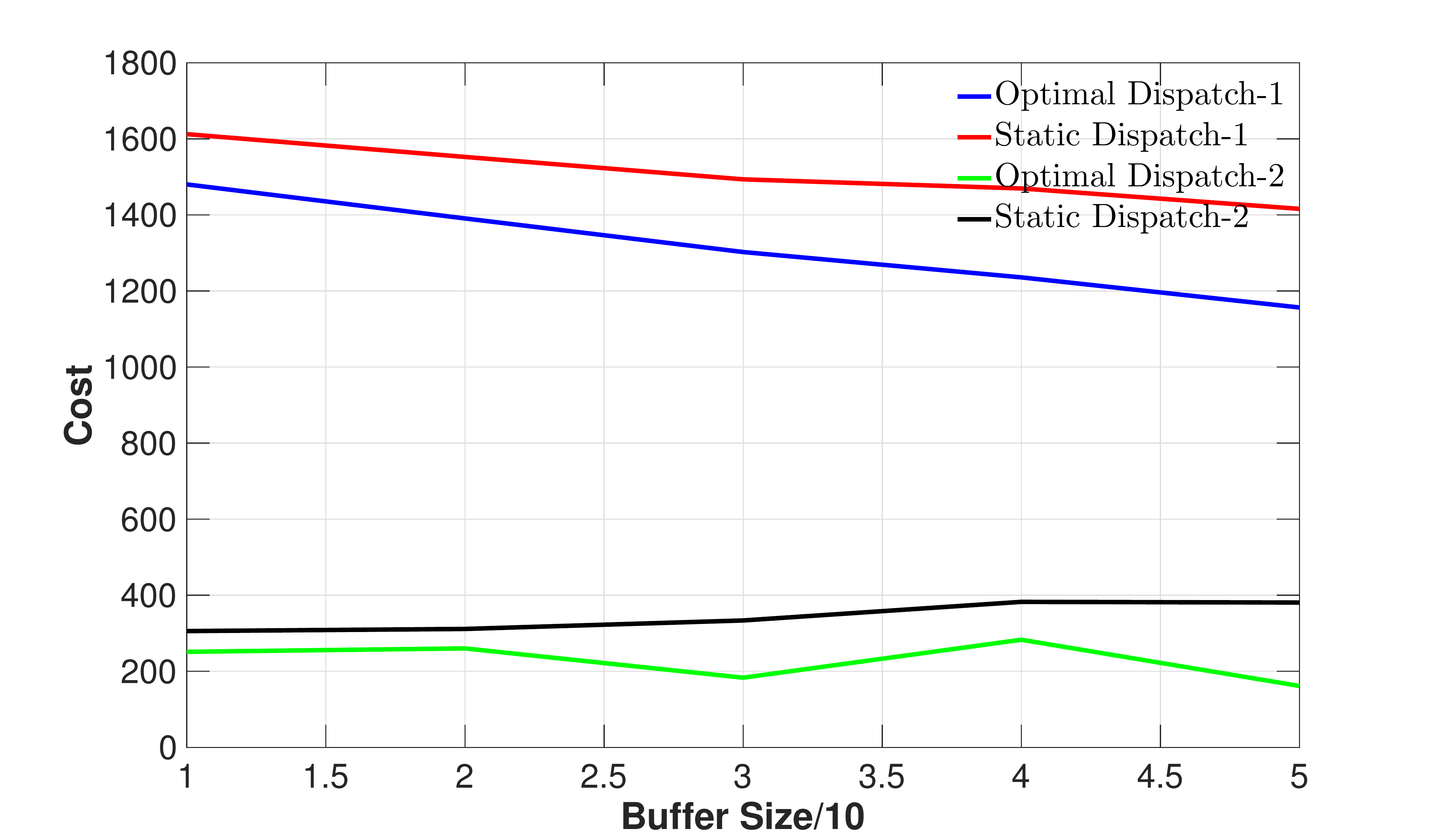}
\vspace{0in}
\caption{Cost as wind availability at $t=1$, and storage buffer at windfarms are increased. Buffer capacity in $i$-th simulation is $10i$, while wind energy $W(1)$ is $i$. Temperature conditions and agents are as in Table~\ref{table1}. 
}
\label{fig10}
\end{figure}

Figures~\ref{fig8}-\ref{fig10} compare the costs averaged over multiple wind realizations of the two policies under various scenarios, for the two schemes. 
Thermal loads are allowed to become energy producers, while wind-farm operators are allowed to store energy in case there is excess energy supply in the market, showcasing how potential prosumer behavior in energy market. 
The particular prices and power generations for Scenario 1, are shown in Table~\ref{table3}.
 

\begin{table}[h]
\vspace{0in}
\includegraphics[width=1\linewidth]{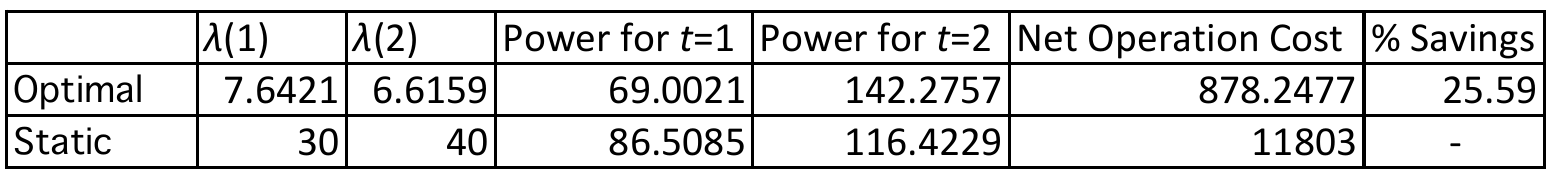}
\caption{Prices, power generation and cost savings.}
\label{table3}
\end{table}
\section{Concluding Remarks} \label{sec-concluding}
The problem of maximizing the social welfare
of a collection of dynamic stochastic agents 
is more complex than stochastic control since agents do
not know the dynamical equations or
utility functions of others. 
It is further complicated by its dynamic, stochastic, decentralized nature,
since each agent's optimal choices depend on the probability distributions
of future prices, which are affected by
the unknown states and actions of all agents.
Yet agents have to make decisions in real-time, as does the ISO since it needs to set prices before agents can decide.

We have exhibited iterative bidding schemes
that attain the optimal performance of a centralized control policy that is aware of the
dynamics, utilities, uncertainties and states of all agents, under appropriate compactness-convexity or LQG assumptions.
It yields the optimal stochastic dynamic locational marginal prices.

The ISO critically exploits the 
sequential information obtained \emph{during} the iterative price-bid process 
to determine the optimal prices and generation/consumption allocations.
This is the stochastic dynamic analog of bidding demand/supply curves
in simple static settings, whence the ISO
can simply intersect cumulative demand and production curves to determine the optimal price.

The social-welfare optimality can potentially result in 
significant economic benefits
in energy markets.
The results may be of interest to general equilibrium theory.

While the agents are all presumed to be ``price takers,"
the scheme can be expected to have some strategic robustness under some monotonicity
assumptions. For example, in the static deterministic case, agents do not benefit from overbidding/underbidding
which drives the price up/down, leading to net losses for the agent in either case.
Examining this in a broader context, while shortening the bid
iteration process, is an important issue.

\section*{Acknowledgment}
The authors thank Pravin Varaiya for identifying a significant error
in an earlier version of the paper.



\bibliographystyle{IEEEtran}

\vspace{-0.25in}
\begin{IEEEbiography}[{\includegraphics[width=1in,height=1.25in,clip,keepaspectratio,angle=270]{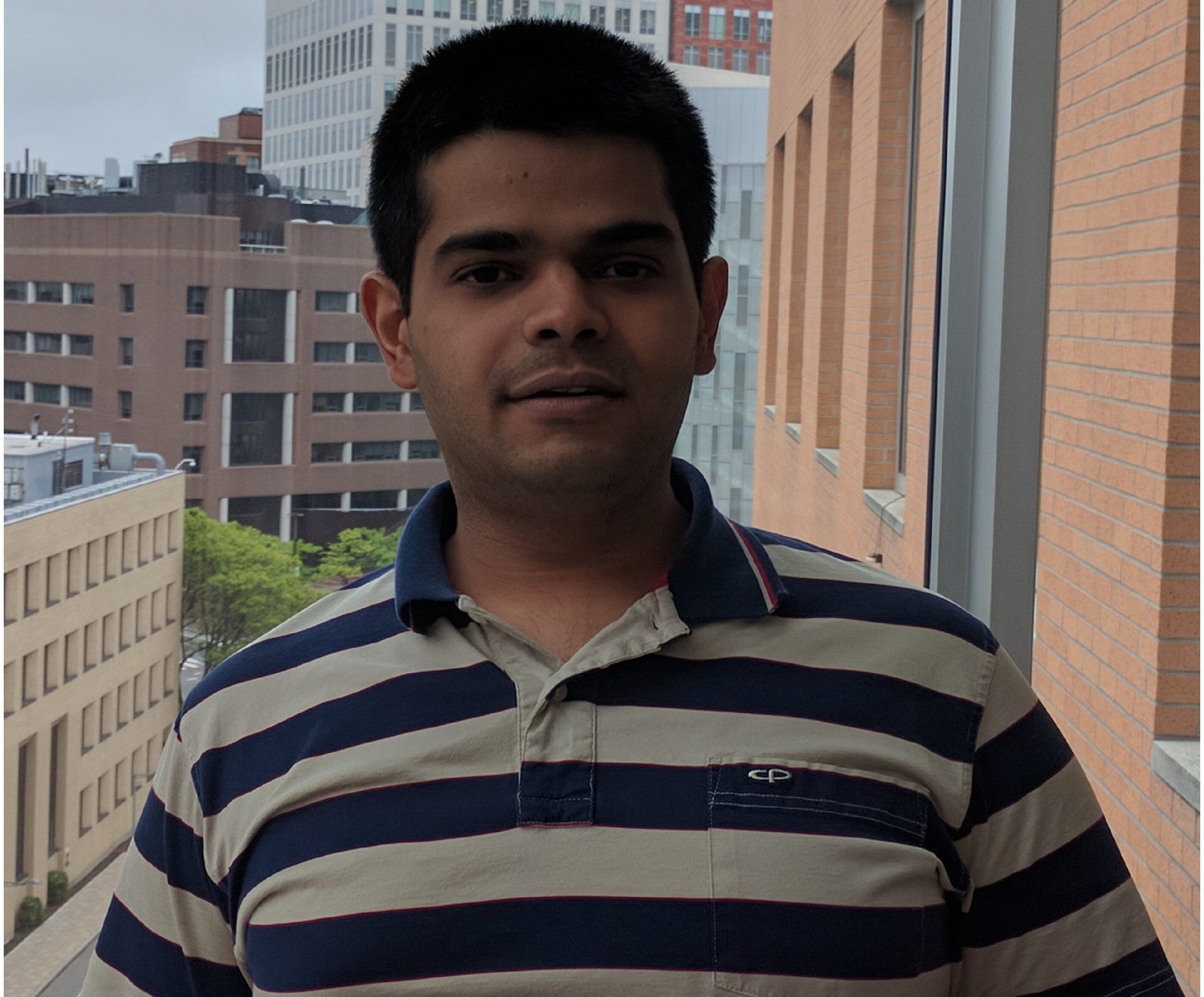}}]{Rahul Singh} received the B.E. degree in electrical engineering from
Indian Institute of Technology, Kanpur, India, in 2009, the M.Sc. degree in Electrical Engineering from University of Notre Dame, South Bend, IN, in 2011, and the Ph.D. degree in electrical and computer engineering from the Department of Electrical and Computer Engineering Texas A\&M University, College Station, TX, in 2015.

He is currently a Postdoctoral Associate at the Laboratory for Information Decision Systems (LIDS), Massachusetts Institute of Technology. His research interests include decentralized control of large-scale complex cyberphysical systems, operation of electricity markets with renewable energy, and scheduling of networks serving real time traffic.
%
\end{IEEEbiography}
\vspace{-0.4in}
\begin{IEEEbiography}[{\includegraphics[width=1in,height=1.25in,clip,keepaspectratio]{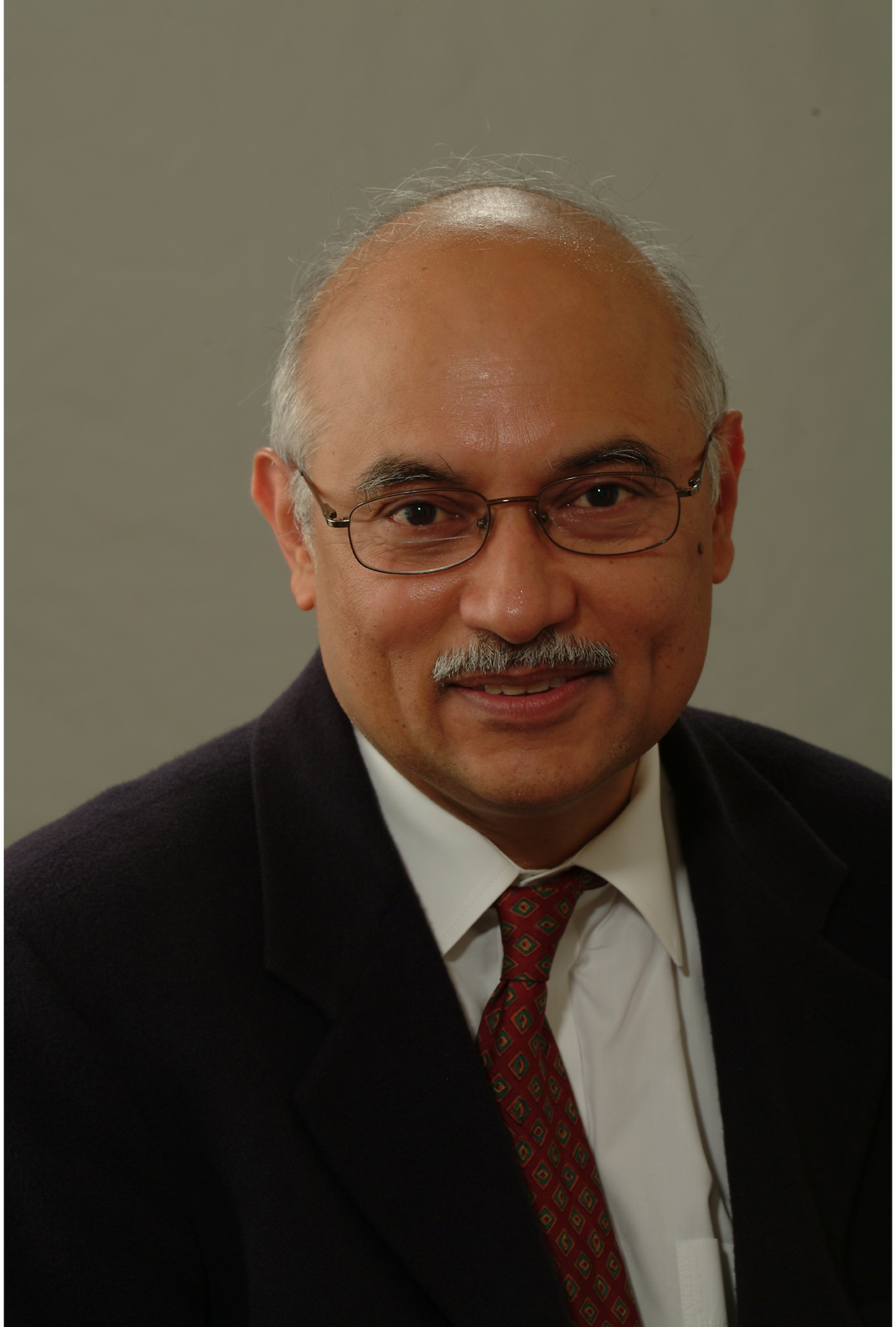}}]{P.~R.~Kumar} B. Tech. (IIT Madras, `73), 
D.Sc. (Washington University, St.~Louis, `77), was a faculty member at UMBC (1977-84) and Univ.~of Illinois, Urbana-Champaign (1985-2011). He is currently at Texas A\&M University. His current research is focused on stochastic systems, energy systems, wireless networks, security, automated transportation, and cyberphysical systems. 

He is a member of the US National Academy of Engineering and The World Academy of Sciences. He was awarded a Doctor Honoris Causa by ETH, Zurich. He as received the
IEEE Field Award for Control Systems, the Donald~P.~Eckman Award of the AACC,  Fred~W.~Ellersick Prize of the IEEE Communications Society, the Outstanding Contribution Award of ACM SIGMOBILE, the Infocom Achievement Award, and the
SIGMOBILE Test-of-Time Paper Award. He is a Fellow of IEEE and ACM Fellow. He was Leader of the Guest Chair Professor Group on Wireless Communication and Networking at Tsinghua University, is a D. J. Gandhi Distinguished Visiting Professor at IIT Bombay, and an Honorary Professor at IIT Hyderabad. He was awarded the Distinguished Alumnus Award from IIT Madras, the Alumni Achievement Award from Washington Univ., and the Daniel Drucker Eminent Faculty Award from the College of Engineering at the Univ.~of Illinois.
\end{IEEEbiography}
\vspace{-0.4in}
\begin{IEEEbiography}[{\includegraphics[width=1in,height=1.25in,clip,keepaspectratio]{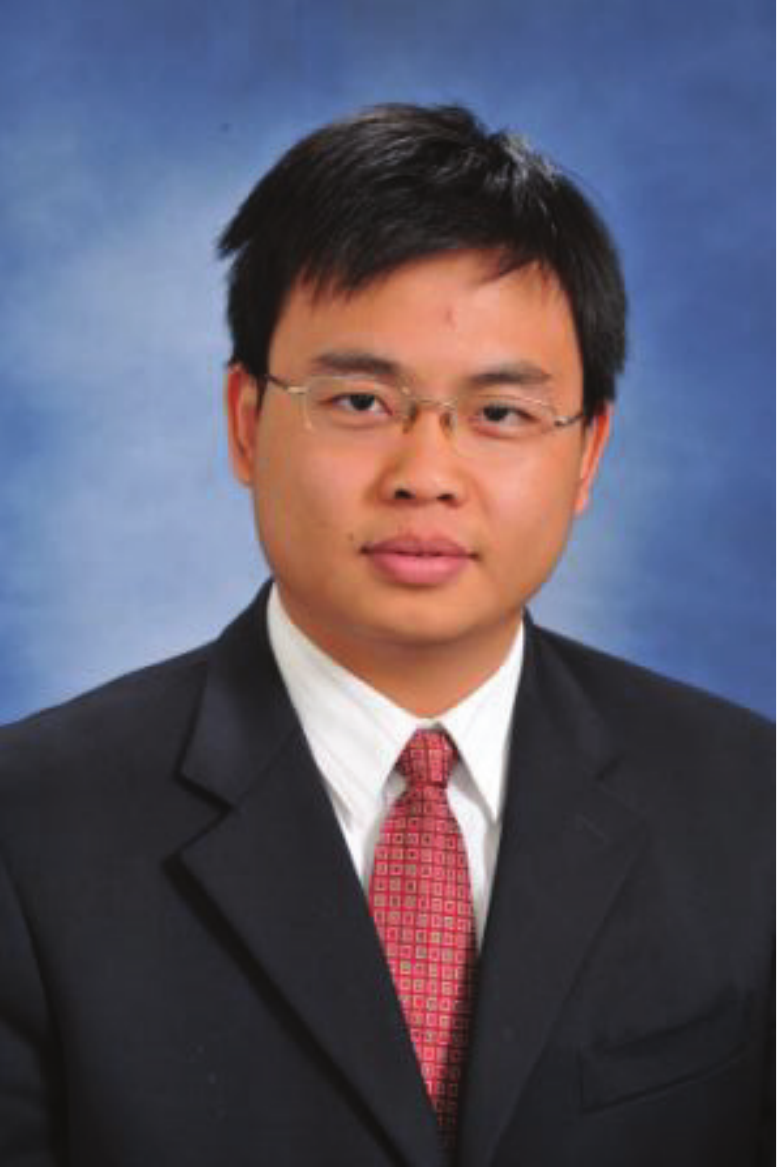}}]{Le Xie} (S'05-M'10-SM'16) received the B.E. degree in electrical engineering from
	Tsinghua University, Beijing, China, in 2004, the M.Sc. degree in engineering
	sciences from Harvard University, Cambridge, MA, in 2005, and the
	Ph.D. degree in electrical and computer engineering from Carnegie
	Mellon University, Pittsburgh, PA, in 2009.
	
	He is currently an Associate Professor with the Department of
	Electrical and Computer Engineering, Texas A\&M University, College
	Station. His
	research interest includes modeling and control of large-scale
	complex systems, smart grid application with renewable energy
	resources, and electricity markets.
	\end{IEEEbiography}




\end{document}